\documentclass[11pt]{article}
\usepackage[utf8]{inputenc}
\usepackage[british]{babel}
\usepackage{cmap}
\usepackage{lmodern}
\usepackage[T1]{fontenc}

\usepackage{amssymb, amsmath, amsthm}
\usepackage[a4paper,margin=1in]{geometry}
\usepackage{ragged2e}

\usepackage{authblk} 
\usepackage{pifont}
\usepackage{graphicx}
\usepackage[dvipsnames,svgnames,table]{xcolor}
\usepackage[figuresright]{rotating}
\usepackage{xtab} 
\usepackage{longtable} 
\usepackage{multirow}
\usepackage{footnote}
\usepackage[stable]{footmisc}
\usepackage{chngpage} 
\usepackage{pdflscape} 
\usepackage[nottoc,notlot,notlof]{tocbibind} 

\usepackage{pgfplots}
\pgfplotsset{every tick label/.append style={font=\footnotesize}}
\pgfplotsset{compat=1.18}
\usepackage{setspace}

\usepackage{array}
\newcolumntype{K}[1]{>{\centering\arraybackslash$}p{#1}<{$}}

\makesavenoteenv{tabular}
\usepackage{tabularx}
\usepackage{booktabs}
\usepackage{threeparttable}
\usepackage[referable]{threeparttablex} 
\newcolumntype{R}{>{\raggedleft\arraybackslash}X}
\newcolumntype{L}{>{\raggedright\arraybackslash}X}
\newcolumntype{C}{>{\centering\arraybackslash}X}
\newcolumntype{A}{>{\columncolor{gray!25}}C}
\newcolumntype{a}{>{\columncolor{gray!25}}c}

\newlength{\tablen}

\usepackage{dcolumn} 
\newcolumntype{.}{D{.}{.}{-1}}

\usepackage{tikz}
\usetikzlibrary{arrows, automata, calc, matrix, patterns, positioning, trees}
\usepackage[semicolon]{natbib} 
\usepackage[hyphens]{xurl}
\usepackage[nopatch=footnote]{microtype}
\usepackage[justification=centerfirst]{caption} 

\usepackage[labelformat=simple]{subcaption}

\DeclareCaptionLabelFormat{parenthesis}{(#2)}
\makeatletter
\renewcommand\p@subfigure{\arabic{figure}.}
\makeatother

\DeclareCaptionLabelFormat{parenthesis}{(#2)}
\makeatletter
\renewcommand\p@subtable{\arabic{table}.}
\makeatother

\usepackage{enumitem}

\setlist[itemize]{leftmargin=2.5\parindent}
\setlist[enumerate]{leftmargin=2.5\parindent}

\usepackage{hyperref} 
\hypersetup{
  colorlinks   = true,    		
  urlcolor     = blue,    		
  linkcolor    = blue,    		
  citecolor    = ForestGreen	
}

\def\addlegendimage{\csname pgfplots@addlegendimage\endcsname}

\theoremstyle{plain}

\newtheorem{corollary}{Corollary}
\newtheorem{lemma}{Lemma}
\newtheorem{proposition}{Proposition}
\newtheorem{theorem}{Theorem}

\theoremstyle{definition}

\newtheorem{definition}{Definition}
\newtheorem{example}{Example}

\theoremstyle{remark}

\makeatletter
\let\@fnsymbol\@alph
\makeatother

\def\keywords{\vspace{.5em} 
{\noindent \textit{Keywords}: }}

\def\AMS{\vspace{.5em} 
{\noindent \textbf{\emph{MSC} class}: }}

\def\JEL{\vspace{.5em} 
{\noindent \textbf{\emph{JEL} classification number}: }}

\newcommand{\HHI}{\textsc{Max-HHI}}

\title{The maximum concentration of results \\ in incomplete round-robin tournaments}

\author{Gergely Cs\'aji\thanks{~Email: \emph{csajigergely@sztaki.hu} \newline Institute for Computer Science and Control (SZTAKI), Hungarian Research Network (HUN-REN), Laboratory on Engineering and Management Intelligence, Budapest, Hungary
\newline  Budapest University of Technology and Economics, Department of Computer Science and Information Theory, Faculty of Electrical Engineering and Informatics, Hungary}
$\qquad \qquad$
\href{https://sites.google.com/view/laszlocsato}{L\'aszl\'o Csat\'o}\thanks{~Corresponding author. Email: \emph{laszlo.csato@sztaki.hun-ren.hu} \newline
Institute for Computer Science and Control (SZTAKI), Hungarian Research Network (HUN-REN), Laboratory on Engineering and Management Intelligence, Research Group of Operations Research and Decision Systems, Budapest, Hungary \newline
Corvinus University of Budapest (BCE), Institute of Operations and Decision Sciences, Department of Operations Research and Actuarial Sciences, Budapest, Hungary}}

\date{\today}

\def\Dedication{
{\noindent
``\emph{The arguments above make it clear that the varying upper and lower bounds of HHI need to be taken into account to obtain a measure of competitive balance that is generally applicable for comparisons over time or across different leagues}'' \citep[p.~293]{OwenRyanWeatherston2007}.
}
\vspace{0.5cm}

\justify }

\begin{document}
\maketitle

\thispagestyle{empty}
\Dedication

\begin{abstract}
\noindent
The Herfindahl--Hirschman Index (HHI), a standard measure of market concentration, is widely used to quantify competitive balance in sports tournaments. To enable comparisons over time or across leagues, both its minimum and maximum should be known. The minimum is usually trivial. Our paper provides results on the maximum HHI in incomplete round-robin tournaments, where the set of matches is represented by a regular but not complete undirected graph.
We prove that the lower bound of maximum HHI is reached when the tournament consists of independent round-robin groups, while the upper bound of maximum HHI is reached when the set of matches is represented by a bipartite graph. Maximising HHI is computationally difficult even for highly structured graph types corresponding to real-world tournaments, where the teams are partitioned into pots and play the same number of matches against opponents from every pot. The complexity of the problem does not depend on the scoring system or the number of opponents from a pot.


\AMS{62P20, 90-08, 91-08}

\JEL{C44, D63, L83, Z20}
\end{abstract}

\clearpage
\restoregeometry

\section{Introduction} \label{Sec1}

Competitive balance, the inequality observed in match outcomes in a tournament, is a fundamental concept of sports economics \citep{Szymanski2003}. Sports have a unique economic structure where, besides competing against each other, multiple teams jointly create another ``product'', the overall league competition. Therefore, in order to keep fans engaged and maximise the entertainment value of the league, it is essential in all leagues to maintain balance and prevent the long-term domination of some teams.

There are two basic tournament formats: the knockout and round-robin competitions.
In a knockout competition, a loser is immediately eliminated; thus, competitive balance is always perfect.
In a round-robin competition, every team plays the same number of matches against any other team, and competitive balance becomes a crucial issue. 

Obviously, competitive balance has many dimensions, which implies the existence of various competitive balance measures \citep{Humphreys2002, FortMaxcy2003, CsatoPetroczy2025a}. Since competitive balance is closely related to the inequality of match outcomes, standard concentration measures from the literature on income distribution and industrial organisation, such as the Gini coefficient \citep{SchmidtBerry2001} and relative entropy \citep{Horowitz1997}, are widely used to measure competitive balance. Another well-known measure from this category is the \emph{Herfindahl--Hirschman index} ($\mathit{HHI}$) \citep{Herfindahl1950, Hirschman1945, Hirschman1964}.
As a measure of competitive balance, $\mathit{HHI}$ is the sum of the squares of teams' shares from the total points allocated in a league \citep{Depken1999}. The number of teams in a round-robin tournament affects both the lower \citep{Depken1999} and upper bounds \citep{OwenRyanWeatherston2007, Avila-CanoRuiz-SepulvedaTriguero-Ruiz2021} of the Herfindahl--Hirschman index; hence, it is usually normalised based on these two extremes to lie in the $\left[ 0,1 \right]$ unit interval \citep{OwenOwen2022, Avila-CanoOwenTriguero-Ruiz2023}.

In recent years, another tournament format, the \emph{incomplete round-robin}, has also become popular in the world of sports. The Union of European Football Associations (UEFA) has changed the design of its club competitions to an incomplete round-robin from the 2024/25 season, and will organise its national team competitions in this format from 2028 \citep{UEFA2026}.

In an incomplete round-robin tournament, each team plays the same number of matches, but they play only against a subset of the other teams \citep{Cochran1971, DevriesereCsatoGoossens2025}. In contrast to the Swiss-system \citep{Csato2017a, SauerCsehLenzner2024}, the schedule of matches is fixed in advance. The incomplete round-robin tournament is effective in reducing travel distance and capacity violations \citep{DevriesereGoossens2025, LiVanBulckGoossens2025}, as well as in avoiding non-competitive matches \citep{CsatoGyimesi2026b, DevriesereGoossensSpieksma2026}.

If the teams play the same number of matches, the minimum $\mathit{HHI}$ is the reciprocal of the number of teams when all matches result in a draw.
However, regarding the other extreme, several interesting questions arise for incomplete round-robin tournaments: What is the maximum value of $\mathit{HHI}$? What is the computational complexity of finding this upper bound? Can reasonable restrictions simplify the problem?
The current paper attempts to answer these questions.

\paragraph{Our contributions}
First, we characterise tournaments with the highest and lowest maximum $\mathit{HHI}$ values. The highest attainable $\mathit{HHI}$ in an incomplete round-robin tournament with $N$ teams, under either the $(2,1,0)$ or $(3,1,0)$ scoring system, is $2/N$, with equality precisely when the set of matches is represented by a bipartite graph (Theorem~\ref{Theo1}). Intuitively, bipartite graphs allow half of the teams to win all their matches, concentrating points as much as possible. At the opposite extreme, the lowest maximum $\mathit{HHI}$ is attained when the tournament consists of disjoint round-robin groups (Theorem~\ref{Theo3}), as such a structure restricts the concentration of wins by confining them to small groups. If the number of teams does not allow creating such equal groups, the lower bound breaks down (Proposition~\ref{PropLowerBoundNondivisibleCounterexample}). Our proof uses equitable graph colouring and majorisation to show that any other match schedule allows an at least as unequal distribution of points. Theorem~\ref{Theo2} also provides a formula for computing this extremal value.

We then investigate the computational complexity of maximising $\mathit{HHI}$ for a given tournament graph. Although the objective involves squared scores, the main difficulty resides in the combinatorial dependencies between match outcomes. We prove $\mathcal{NP}$-hardness even for highly structured graph types corresponding to real-world tournaments, where the teams are partitioned into pots and play the same number of matches against opponents from every pot. In particular, Theorem~\ref{Theo5} shows this for the $(\omega,1,0)$ scoring system with $\omega \geq 3$, while Theorem~\ref{Theo7} proves hardness for the $(2,1,0)$ scoring system.
If the total number of points is fixed (draws are prohibited or the $(2,1,0)$ scoring system is applied), maximising $\mathit{HHI}$ amounts to maximising the sum of squared wins, which we relate to the minimum odd-cycle transversal problem through a reduction from Not-All-Equal 3-SAT. The $(3,1,0)$ scoring system is more challenging to analyse, since draws change both the numerator and denominator of $\mathit{HHI}$. Our reduction encodes a satisfiability instance using variable cycles and clause triangles: alternating winning patterns represent truth assignments, and clause gadgets enforce satisfaction. Finally, carefully chosen auxiliary graphs are used to translate this combinatorial structure into the fractional $\mathit{HHI}$ objective. 


\bigskip
Our study is not the first to examine incomplete round-robin tournaments from the perspective of graph theory. \citet{LiVanBulckGoossens2025} provide theoretical results on the complexity of constructing an incomplete round-robin tournament. \citet{deWerraUrrutiaAssuncao2026} consider the problem of minimising the number of breaks (two consecutive home or two consecutive away games); \citet{DevriesereGoossensSpieksma2027} aim to find break-free timetables that minimise venue capacity violations. \citet{WesselSpieksma2026} establish the complexity of determining whether there exists a fair draw where the strengths of the opponent sets are equal for each team.
However, none of these previous works involves finding match outcomes to optimise an objective function, which we will do. 

The paper is organised as follows.
Section~\ref{Sec2} introduces the problem. Section~\ref{Sec3} derives strict lower and upper bounds on the maximum $\mathit{HHI}$ in incomplete round-robin tournaments. Section~\ref{Sec4} deals with the complexity of computing the maximum $\mathit{HHI}$.
Finally, Section~\ref{Sec5} concludes.

\section{Preliminaries} \label{Sec2}



Let $G = (V,E)$ be a $d$-regular undirected graph representing the set of matches played in an incomplete round-robin tournament, where the vertices are the teams and the edges are the matches.
Each match between teams $v$ and $u$ has three possible outcomes: team $v$ wins, team $v$ loses, or the result is a draw.
Denote the number of teams by $\lvert V \rvert = N$.

For each vertex $v \in V$, let $w_v, d_v, \ell_v$ denote the number of wins, draws, and losses, respectively, such that $w_v + d_v + \ell_v = \text{deg}(v)$. If scoring system $(\omega,1,0)$ is used to rank the teams, the score of vertex $v$ is given by $p_v = \omega w_v + d_v$, where $\omega$ is the value of a win.

We aim to find a valid orientation of outcomes on $E$ that maximises competitive imbalance, measured by the Herfindahl--Hirschman index ($\mathit{HHI}_\omega$) of the tournament. $\mathit{HHI}_\omega$ is the sum of squared scores divided by the square of total points:
\[
\mathit{HHI}_\omega(G) = \frac{\sum_{v \in V} p_v^2}{P^2},
\]
where $P = \sum_{v \in V} p_v$ is the total number of points awarded in the tournament.

The lower bound of $\mathit{HHI}_\omega$, associated with perfect competitive balance, is $1/N$ in any graph $G$ and for any scoring system $(\omega,1,0)$. It can be easily achieved if all teams collect the same number of points by playing only draws.
The upper bound of $\mathit{HHI}_\omega$ is denoted by $\mathit{HHI}_\omega^\ast$.

Let \HHI\ denote the optimisation problem that aims to find the maximum attainable $\mathit{HHI}_\omega$ of a given input graph $G$ and scoring system $(\omega,1,0)$.

Regularity is a natural requirement for a tournament graph, since ranking would be challenging if the number of matches played by the teams differs \citep{LeivaBertran2025}. However, the organisers usually impose other constraints in order to avoid high differences between the strengths of opponents for two teams. This can be achieved by creating pots based on the strengths of the teams, and drawing the matches such that each team plays the same number of matches against teams from each pot, possibly without playing against teams from its own pot.

\begin{figure}[t!]
\centering
\begin{tikzpicture}[scale=1, auto=center, transform shape, >=triangle 45]
\tikzstyle{every node}=[draw,shape=circle];
  \node (n1) at (120:2) {};
  \node (n2) at (60:2) {};
  \node (n3) at (180:2) {};
  \node (n4) at (0:2) {};
  \node (n5) at (240:2) {};
  \node (n6) at (300:2) {};

  \foreach \from/\to in {n1/n2,n1/n4,n1/n6,n2/n3,n2/n5,n3/n4,n3/n6,n4/n5,n5/n6}
    \draw (\from) -- (\to);
\end{tikzpicture}

\caption{The graph corresponding to the group stage of the Romanian football cup}
\label{Fig1}
\end{figure}
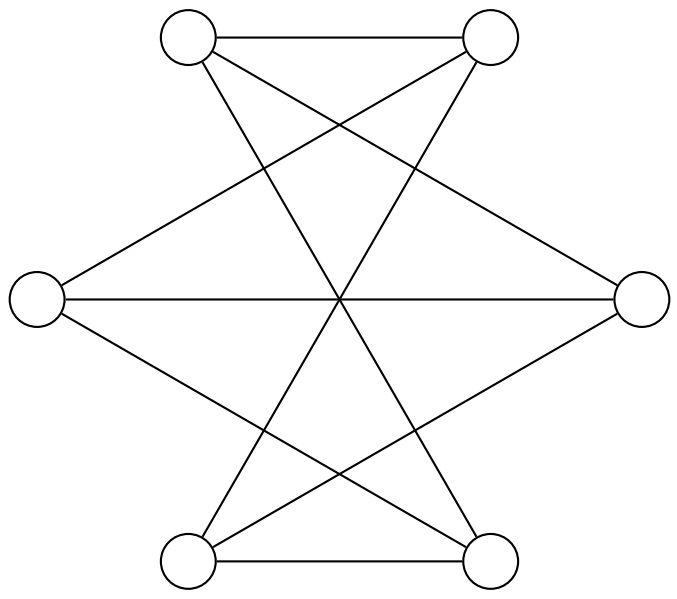

For example, the Romanian football cup (Cupa Rom{\^ a}niei) contains a ``group stage'' since the 2022/23 season. Each of the four groups contains two teams from Pots 1--3, respectively. Every team plays one match against one team from each of the three pots in their group.
Figure~\ref{Fig1} shows the graph associated with a group in this tournament, which is unique if the identities of the teams are ignored.
The maximum of the Herfindahl--Hirschman index is reached under the point distribution (9,9,9,0,0,0), when half of the teams win all their matches.

Let $G=(V,E)$ and $V=V_1\mathbin{\dot\cup}\cdots\mathbin{\dot\cup}V_k$ be a partition into classes of equal size. Denote the neighbours of $v$ by $N(v)$ for all $v \in V$.
\begin{enumerate}
    \item The graph (and the partition) is \emph{fully balanced with parameter $r$} if for all $v \in V$ and $1 \leq j \leq k$:
    \[
    |N(v)\cap V_j| = r
    \]
    \item The graph (and the partition) is \emph{externally balanced with parameter $r$} if for all $v \in V_i$ and $j \neq i$:
    \[
    |N(v)\cap V_j|=r
    \]
\end{enumerate}

Fully balanced graphs represent tournaments where each team plays exactly $r$ matches with teams from every pot, including its own. Externally balanced graphs differ only in that teams do not play against other teams from their own pot, but play exactly $r$ matches with teams from every other pot.

Both structures arise naturally in major European sports competitions since the 2023/24 season.
In the league phase of the UEFA Champions League and Europa League, the 36 participating clubs are divided into four seeding pots of nine clubs each, and every club plays against two opponents from each pot, including its own. Thus, the set of matches yields a fully balanced graph with $r=2$. The UEFA Conference League employs six pots of six clubs, with each club facing exactly one opponent from each pot, giving a fully balanced graph with $r=1$.
In the Champions Hockey League, the 24 participating teams are divided into four seeding pots of six teams each, and every team faces exactly two opponents from each of the other three pots. Thus, the set of matches is represented by an externally balanced graph with parameter $r=2$.

\section{Bounds on the maximum Herfindahl--Hirschman index} \label{Sec3}

In this section, we study $\mathit{HHI}$ in incomplete round-robin tournaments, and show that it has a strict, structurally bounded global maximum. In particular, Section~\ref{Sec31} establishes an upper bound on $\mathit{HHI}_\omega^\ast$, Section~\ref{Sec32} computes $\mathit{HHI}_\omega^\ast$ in disjoint round-robin tournaments, and Section~\ref{Sec33} establishes a lower bound on $\mathit{HHI}_\omega^\ast$.
Recent results about the maximum $\mathit{HHI}$ in round-robin tournaments are summarised in Appendix~\ref{Appendix_A}.

Throughout this section, $\omega\geq2$ is assumed to be an exogenous integer. The case $\omega=3$, which corresponds to the standard football scoring system, is stated explicitly whenever the general formula is less transparent.

\subsection{Upper bound and characterisation} \label{Sec31}

The following result provides a strict upper bound on the maximum Herfindahl--Hirschman index for incomplete round-robin tournaments.

\begin{theorem} \label{Theo1}
Let $G$ be a $d$-regular graph representing an incomplete round-robin tournament. Then
\[
\mathit{HHI}_\omega^\ast(G)\leq \frac{C_\omega}{N},
\]
where
\[
C_\omega=
\begin{cases}
2, & \omega\in\{2,3\},\\[1mm]
\displaystyle\frac{(\omega^2-2)^2}
{4\omega(\omega-1)(\omega-2)}, & \omega\geq4.
\end{cases}
\]
For $\omega\in\{2,3\}$, equality holds if and only if $G$ is bipartite. For
$\omega\geq4$, put
\[
\rho_\omega:=
\frac{2(\omega^2-2\omega+2)}{(\omega^2-2)(\omega-2)}.
\]
Equality holds if and only if $G$ is the disjoint union of a $d$-regular bipartite graph on $\rho_\omega N$ vertices and a $d$-regular graph on the remaining vertices; no match in the first graph is a draw and all matches in the second graph are draws.
\end{theorem}

\begin{proof}
Let $W$ be the total number of wins in the tournament. Because $G$ is $d$-regular, the number of edges is $Nd/2$. The total sum of draw-degrees is $Nd-2W$. Thus, the sum of points is
\begin{equation}
\label{eq:numpoints}
P=\omega W+(Nd-2W)=Nd+(\omega-2)W.    
\end{equation}

Using
$w_v+d_v+\ell_v=d$, we get $w_v^2 = d w_v - w_v d_v - w_v \ell_v$ and $d_v^2 = d d_v - w_v d_v - \ell_v d_v$.
Expanding the squared points  $(\omega w_v+d_v)^2$ and using $\sum_v w_v=W$ and $\sum_v d_v=Nd-2W:$
\begin{equation}
\label{eq:numsquares}
S=Nd^2+d(\omega^2-2)W-\Gamma_\omega,
\end{equation}
where
\[
\Gamma_\omega:=\sum_{v\in V}
\left(\omega^2w_v\ell_v+(\omega-1)^2w_vd_v+d_v\ell_v\right)\geq0.
\]

Put $x:=W/(Nd)$. Since $0\leq W\leq Nd/2$, we have
$0\leq x\leq1/2$, and consequently
\begin{equation}
\mathit{HHI}_\omega(G)
\leq
\frac{1}{N}
\frac{1+(\omega^2-2)x}{(1+(\omega-2)x)^2}.
\label{EqGeneralUpperEnvelope}
\end{equation}

The derivative of this fraction as a function of $x$ has the same sign as
\[
\omega^2-2\omega+2-(\omega^2-2)(\omega-2)x.
\]
For $\omega\in\{2,3\}$, the expression in
\eqref{EqGeneralUpperEnvelope} is increasing on $[0,1/2]$, and its value at
$x=1/2$ is $2/N$. In particular, for $\omega \in \{ 2,3\}$ we get
\[
\mathit{HHI}_\omega^\ast(G)\leq\frac{2}{N}.
\]
For $\omega\geq4$, the unique maximum is attained at
\[
x_\omega^\ast=
\frac{\omega^2-2\omega+2}{(\omega^2-2)(\omega-2)},
\]
and substituting this value into \eqref{EqGeneralUpperEnvelope} gives
$C_\omega/N$.

Equality in the envelope requires $\Gamma_\omega=0$. Hence every vertex is a pure winner, a pure drawer, or a pure loser. For $\omega\in\{2,3\}$, equality also requires $x=1/2$, so there are no draws; the pure winners and pure losers form the two sides of a bipartition. Conversely, a pure bipartite orientation attains $2/N$.

For $\omega\geq4$, equality requires $W=x_\omega^\ast Nd$. The pure winners and pure losers therefore form a $d$-regular bipartite union on $2W/d=\rho_\omega N$ vertices, while the pure drawers induce a $d$-regular union on the remaining vertices. There can be no edges between these two parts. Conversely, orienting the first part bipartitely and making every match a draw in the second part attains the maximum.
\end{proof}

\begin{corollary}
Letting each match be a draw gives a $C_\omega$-approximation to \HHI, since this configuration has $\mathit{HHI}_\omega=1/N$. In particular, it is a 2-approximation for $\omega\in\{2,3\}$.
\end{corollary}

Theorem~\ref{Theo1} is not only a theoretical curiosity, as there exist incomplete round-robin tournaments where the matches are given by a bipartite graph. We have already discussed the Romanian football cup in Section~\ref{Sec2}.
The Leagues Cup is an international football competition for Major League Soccer (MLS) and Mexican clubs. It starts with a league phase since 2025, contested by 18-18 MLS and Mexican clubs, which are divided into six groups of three teams each, determined by three tiers based on historical performances and two geographical regions. In the groups, each team plays against the three teams from the other country; thus, each match is played by one club from MLS and one from Mexico. 

\subsection{The maximum $\mathit{HHI}$ in disjoint round-robin tournaments} \label{Sec32}

In this section, we consider $\bigcup K_{d+1}$. For $\omega=3$, the relevant component states are precisely the truncated cascade distributions studied by \citet{Avila-CanoRuiz-SepulvedaTriguero-Ruiz2021}, and their $\mathit{HHI}$ is given by the formula stated in Section~\ref{Sec2}. The same score pattern extends to general $\omega$. Since different components of a disjoint union may be in two adjacent states, we need their total points and sums of squares separately.

Starting from the all-draw configuration, an \emph{upgrade} of a vertex means replacing its draws against the current bottom tier of its component by wins, one after another. A \emph{greedy win allocation} always chooses a vertex having the minimum possible number of losses, completes its upgrade, and then continues with another such vertex.

The following lemma gives the complete structural reduction that we need.

\begin{lemma}[Structure of an optimal configuration] \label{lem:complete-helper}
There exists a configuration maximising the $\mathit{HHI}$ on $G_{\mathrm{worst}}=\bigcup K_{d+1}$ that is a greedy win allocation in which every started upgrade is completed. Consequently, every component is in a truncated cascade state: it has singleton upper tiers and one bottom draw tier. Moreover, the number of fully upgraded vertices in any two components differs by at most one.
\end{lemma}

\begin{proof}
For $\omega=2$, the total number of points is fixed. If a highest-scoring vertex draws against another vertex, changing that draw into a win of the highest-scoring vertex transfers one point from the latter vertex to the former and strictly increases the sum of squared scores. Repeating this argument for second highest scorers, etc., gives the complete cascade in every component. Hence, all assertions are immediate for $\omega=2$.

Assume from now on that $\omega\geq3$, and put
\[
c_\omega:=\omega-2.
\]
We use the upper-bound constant $C_\omega$ from
Theorem~\ref{Theo1}. In particular, every configuration on the whole
$N$-vertex graph satisfies
\begin{equation}
\mathit{HHI}_\omega\leq
\min\left\{1,\frac{C_\omega}{N}\right\}.
\label{EqCompleteGlobalUpperBound}
\end{equation}

We first prove a saturation claim. Divide the vertices into tiers based on their scores: $u,v$ belong to the same tier if and only if $p_u=p_v$. Fix a component and suppose that some
vertices of this component (possibly 0) have already been fixed as singleton upper tiers.
Let $R$ be the set of remaining vertices and write
\[
|R|=m+1.
\]
Among the $\mathit{HHI}$-maximising configurations respecting the already fixed upper
tiers, choose one that minimises the number of decisive matches inside $R$.
Let $v$ be a highest-scoring vertex of $R$.

\medskip
\noindent
\textbf{Step 1: The highest remaining vertex has no loss.}

Suppose that $v$ loses to some $u\in R$. Reverse this result, so that $v$
defeats $u$. The total number of points remains unchanged, while the score of
$v$ increases by $\omega$ and the score of $u$ decreases by $\omega$.
Since $p_v\geq p_u$, the corresponding change in the sum of squared scores is
\begin{align*}
&(p_v+\omega)^2+(p_u-\omega)^2-p_v^2-p_u^2=
2\omega(p_v-p_u)+2\omega^2
>0.
\end{align*}
This contradicts optimality. Hence, $v$ has only wins and draws against the
other $m$ vertices of $R$.

Temporarily replace all these matches by draws. Let
\[
r_1\leq r_2\leq\cdots\leq r_m
\]
be the resulting scores of the other vertices. If $v$ defeats a vertex with
baseline score $r_j$ but draws with a vertex having baseline score
$r_i<r_j$, exchange these two results. The score of $v$ and the total points
remain unchanged, while the change in the sum of squared scores is
\begin{align*}
&(r_i-1)^2+r_j^2-r_i^2-(r_j-1)^2=2(r_j-r_i)>0.
\end{align*}
Therefore, if $v$ has exactly $w$ wins, we may assume that it defeats the
vertices with baseline scores $r_1,\ldots,r_w$.

\medskip
\noindent
\textbf{Step 2: A highest remaining vertex either draws every match or wins
every match.}

The proof of Step 2 is intuitive but technical; thus it is provided in Appendix~\ref{Appendix_B}.

\medskip
\noindent
\textbf{Step 3: Recursive construction of a truncated cascade.}

Recall that, subject to the already fixed singleton upper tiers, we chose an $\mathit{HHI}$-maximising configuration having the minimum number of decisive matches inside $R$. If $0<w<m$, the all-draw boundary has fewer decisive matches. Consequently, it cannot have the same $\mathit{HHI}$, by the tie-breaking choice. The all-win boundary must therefore attain the same maximum $\mathit{HHI}$. Hence, there exists an optimal configuration in which $v$ defeats every other vertex of $R$. We fix $v$ as the next singleton upper tier and repeat the argument on
$R\setminus\{v\}$.

The recursion may instead reach a stage at which a highest-scoring remaining
vertex $v$ already has $w=0$. Then $v$ has score $m$. Every other remaining
vertex has score at most $m$, and, since it draws with $v$, its score from
matches inside $R\setminus\{v\}$ is at most $m-1$. On the other hand, the
$m$ vertices of $R\setminus\{v\}$ play $\binom{m}{2}$ matches, each awarding
at least two points, so their total score from these matches is at least
$m(m-1)$. Hence, every one of them has score exactly $m-1$ from the induced
subtournament, and every match in that subtournament awards exactly two
points. Since $\omega\geq3$, all those matches are draws. Thus, the whole
remaining set is one bottom draw tier.

Applying this argument independently to every component gives a $\mathit{HHI}$-maximising configuration in which every component is a truncated cascade.

\medskip
\noindent
\textbf{Step 4: The cascade depths of different components differ by at
most one.}

Suppose, to the contrary, that two components have bottom draw tiers of
sizes $x\ge 1$ and $y$, where$
y\geq x+2.$ 
In the component with bottom size $x$, undo the last completed upgrade. This
removes the $x$ wins by which one vertex was separated from a bottom tier of
size $x+1$. In the component with bottom size $y$, use these $x$ wins to
start upgrading one vertex against exactly $x$ of its $y-1$ bottom-tier
opponents. The total number of decisive matches, and hence the total number
of points, remains unchanged.

For a bottom tier of size $z$, giving one vertex $w$ wins against its peers
changes the sum of squared scores by
\begin{align}
\Delta S(z,w)
={}&
\left[z-1+(\omega-1)w\right]^2-(z-1)^2
\nonumber+
w\left[(z-2)^2-(z-1)^2\right]
\nonumber
\\ =&
\left[2c_\omega(z-1)+1\right]w
+(\omega-1)^2w^2.
\label{EqCompletePartialLayer}
\end{align}
Therefore, the loss caused by undoing the completed upgrade is
$\Delta S(x+1,x)$, while the gain in the other component is
$\Delta S(y,x)$. Their difference is
\begin{align*}
\Delta S(y,x)-\Delta S(x+1,x)
&=
2c_\omega x(y-x-1)>0.
\end{align*}
We have obtained a configuration with the same total points and a strictly larger sum of squared scores, contradicting optimality. Hence, the bottom tier sizes of any two components differ by at most one; equivalently, their number of fully upgraded vertices differs by at most one.
\end{proof}

We are ready to give a formula for the maximum $\mathit{HHI}$ of disjoint round-robin tournaments.

\begin{theorem} \label{Theo2}
Let
$G_{\mathrm{worst}}=\bigcup_{i=1}^{N/(d+1)}K_{d+1}$ and put
$c=N/(d+1)$. Then
$$
\mathit{HHI}_\omega^*(G_{\mathrm{worst}})
=
\frac{d+1}{N}
\max_{x\in\{1,\ldots,d+1\}}
\max_{\alpha\in\{0,1/c,\ldots,1\}}
\frac{\alpha S_c^\omega(x-1)+(1-\alpha)S_c^\omega(x)}
{[\alpha P_c^\omega(x-1)+(1-\alpha)P_c^\omega(x)]^2},
$$
where
\begin{align}
P_c^\omega(x)
&=
\frac{\omega d(d+1)-(\omega-2)x(x-1)}{2},
\label{EqPcGeneralOmega}\\
S_c^\omega(x)
&=
\omega^2\left[
\frac{d(d+1)(2d+1)}{6}
-
\frac{(x-1)x(2x-1)}{6}
\right]
+x(x-1)^2.
\label{EqScGeneralOmega}
\end{align}
Here $x$ is the baseline number of bottom vertices in a component and
$\alpha$ is the fraction of components upgraded to $x-1$ bottom vertices.
\end{theorem}

\begin{proof}
By Lemma~\ref{lem:complete-helper}, an optimal configuration consists of
truncated cascade components whose numbers of bottom vertices differ by at
most one. Let $k$ components have $x-1$ bottom vertices and the remaining
$c-k$ components have $x$ bottom vertices; put $\alpha=k/c$.

A component with $x$ bottom vertices has all matches decisive except for the
$\binom{x}{2}$ draws among the bottom vertices. Since a decisive match awards
$\omega$ points in total and a draw awards two points, its total number of
points $P_c^\omega(x)
$ is
\begin{align*}
\omega\left[
\binom{d+1}{2}-\binom{x}{2}
\right]
+2\binom{x}{2}=
\frac{\omega d(d+1)}{2}
-\frac{\omega x(x-1)}{2}
+x(x-1)=
\frac{\omega d(d+1)-(\omega-2)x(x-1)}{2},
\end{align*}
which gives \eqref{EqPcGeneralOmega}.

Its top vertices have scores
$\omega x,\omega(x+1),\ldots,\omega d$, while every bottom vertex has score
$x-1$. Using
$$
\sum_{j=1}^n j^2=\frac{n(n+1)(2n+1)}{6},
$$
we obtain
\begin{align*}
S_c^\omega(x)
&=
\sum_{j=x}^d(\omega j)^2
+
\sum_{i=1}^x(x-1)^2=
\omega^2\sum_{j=x}^d j^2+x(x-1)^2=
\omega^2\left(
\sum_{j=1}^d j^2-
\sum_{j=1}^{x-1}j^2
\right)
+x(x-1)^2\\
&=
\omega^2\left[
\frac{d(d+1)(2d+1)}{6}
-
\frac{(x-1)x(2x-1)}{6}
\right]
+x(x-1)^2,
\end{align*}
which gives \eqref{EqScGeneralOmega}. This is exactly the truncated cascade
distribution of
\citet{Avila-CanoRuiz-SepulvedaTriguero-Ruiz2021}.

The total sum of squares and total points are
$$
S(\alpha)
=
c\left[
\alpha S_c^\omega(x-1)+(1-\alpha)S_c^\omega(x)
\right], \quad 
P(\alpha)
=
c\left[
\alpha P_c^\omega(x-1)+(1-\alpha)P_c^\omega(x)
\right].
$$
Substitution into $\mathit{HHI}_\omega=S/P^2$ proves the formula. Since the objective
is a ratio of a linear function and the square of a linear function in
$\alpha$, its continuous maximum may lie in the interior; the discrete
maximum is obtained by checking all displayed values of $x$ and $\alpha$.
\end{proof}

\subsection{Lower bound and characterisation} \label{Sec33}

The following result provides a strict lower bound on the maximum Herfindahl--Hirschman index for incomplete round-robin tournaments by the $\mathit{HHI}_\omega^\ast$ of independent round-robin tournaments, computed in Theorem~\ref{Theo2}.

\begin{theorem}\label{Theo3}
Let $G$ be a simple $d$-regular graph on $N$ vertices representing an
incomplete round-robin tournament, where $d\geq 1$ and $d+1$ divides $N$.
Put
\[
c:=\frac{N}{d+1},\qquad K:=K_{d+1},\qquad
G_{\mathrm{worst}}:=\bigcup_{i=1}^{c}K_{d+1}.
\]
Then, we have that
\begin{equation}
\mathit{HHI}_\omega^\ast(G)
\geq \mathit{HHI}_\omega^\ast(G_{\mathrm{worst}})
\geq \frac{d+1}{N}\mathit{HHI}_\omega^\ast(K).
\label{EqLowerBoundCompleteComponents}
\end{equation}
The first inequality is an equality for $\omega\in\{2,3,4,5\}$
if and only if $G\cong G_{\mathrm{worst}}$.
For larger $\omega$, equality may also hold for other graphs.
\end{theorem}

\begin{proof}
Put $c=N/(d+1)$. We first choose a reference configuration on
$G_{\mathrm{worst}}$ and construct a configuration of $G$ whose score
vector majorises the reference score vector. We then prove strictness
for $\omega\in\{2,3,4,5\}$ when $G\not\cong G_{\mathrm{worst}}$.

Choose an optimal configuration of the entire graph
$G_{\mathrm{worst}}$. By Lemma~\ref{lem:complete-helper}, the
reference may be chosen so that every component has $r$ fully upgraded
winners, while $t$ components have one additional fully upgraded winner,
for some $0\leq r\leq d$ and $0\leq t<c$. We take $t=0$ if $r=d$,
since no further nontrivial upgrade is possible.

Let
\(
x:=d+1-r.
\)
The decreasing reference score vector $\mathbf p=(p_1,\ldots,p_N)$
consists of $c$ copies of each of
\[
\omega d,\ \omega(d-1),\ \ldots,\ \omega(d-r+1),
\]
followed, when $t>0$, by $t$ scores equal to $\omega(d-r)$.
All remaining scores are $x-1$ or $x-2$ (only $x-1$ when $t=0$).
The total number of decisive matches in this configuration is
\begin{equation}
W=c\left(rd-\binom r2\right)+t(d-r).
\label{EqReferenceWinsUnion}
\end{equation}

\medskip
\noindent\textbf{Step 1: Equitable colouring.}
By the Hajnal--Szemer\'edi theorem
\citep{HajnalSzemeredi1970,KiersteadKostochka2008},
$G$ has a proper equitable $(d+1)$-colouring, i.e. a colouring with independent colour
classes $A_0,\ldots,A_d$ of size $c$. Order these classes greedily:
choose $A_0$ arbitrarily and, after choosing the first $j$ classes,
select $A_j$ among the remaining classes such that it minimises its number
of neighbours in their union. Write
\[
U_j:=A_0\cup\cdots\cup A_{j-1}\qquad (U_0=\varnothing).
\]

\medskip
\noindent\textbf{Step 2: A bound on the internal edges.}
We claim that
\begin{equation}
e_G(U_j)\leq c\binom j2
\qquad (0\leq j\leq d+1).
\label{EqPrefixInternalEdgeBound}
\end{equation}
It is convenient to measure the difference by
\[
\delta_j:=c\binom j2-e_G(U_j).
\]
Thus $\delta_1=0$. Suppose $\delta_j\geq0$ and $1\leq j\leq d-1$.
Since $G$ is $d$-regular, the number of edges from $U_j$ to the
remaining $d+1-j$ colour classes is
\[
jcd-2e_G(U_j)=jc(d+1-j)+2\delta_j.
\]
The greedy choice of $A_j$ therefore gives
\[
e_G(A_j,U_j)\leq jc+\frac{2\delta_j}{d+1-j}.
\]
Consequently,
\begin{equation}
\delta_{j+1}
\geq\left(1-\frac{2}{d+1-j}\right)\delta_j\geq0.
\label{EqDeficitInduction}
\end{equation}
This proves \eqref{EqPrefixInternalEdgeBound} for $j\leq d$.
The same bound at $j=d+1$ is an equality, because $G$ has
$N d/2=c\binom{d+1}{2}$ edges.

If $t>0$, one additional, partial class is needed. Choose
$B\subseteq A_r$ with $|B|=t$, consisting of the $t$ vertices with
fewest neighbours in $U_r$, and put $U:=U_r\cup B$.
Since $A_r$ is independent, averaging the degrees into $U_r$ gives
\[
\begin{aligned}
e_G(U)
&\leq e_G(U_r)+\frac{t}{c}\,e_G(A_r,U_r)\leq c\binom r2-\delta_r
       +\frac{t}{c}\left(rc+\frac{2\delta_r}{d+1-r}\right)\\
&=c\binom r2+rt-
  \left(1-\frac{2t}{c(d+1-r)}\right)\delta_r\leq c\binom r2+rt.
\end{aligned}
\label{EqPartialClassInternalEdgeBound}
\]
The last step follows from $t<c$ and $d+1-r\geq2$ whenever $t>0$.

\medskip
\noindent\textbf{Step 3: Constructing a majorising score vector.}
For an edge between $A_i$ and $A_j$, $i<j$, call $i$ its
\emph{leading index}. List the edges by nondecreasing leading index.
If $t>0$, list edges led by vertices of $B$ first among those whose
leading index is $r$. Select the first $W$ edges, make their
lower-indexed endpoint the winner, and declare all other matches draws.
Let $\mathbf q=(q_1,\ldots,q_N)$ be the resulting decreasing score vector.

Fix $1\leq j\leq r$. The number of edges incident with $U_j$ is
\[
F_j=jcd-e_G(U_j)
\geq c\left(jd-\binom j2\right)=:W_j.
\]
All these edges precede edges with leading index at least $j$.
Since $W\geq W_j$, at least $W_j$ selected edges have their winner
in $U_j$. At most $e_G(U_j)\leq c\binom j2$ selected edges have
both endpoints in $U_j$.

Initially, if all matches were draws, the teams in $U_j$ would have
$jcd$ points. Changing a crossing draw into a win by a team in $U_j$
adds $\omega-1$ points to $U_j$, whereas changing an internal draw
adds $\omega-2$ points. Hence
\[
\begin{aligned}
\sum_{v\in U_j}q_v
&\geq jcd+(\omega-1)W_j-e_G(U_j)\geq jcd+(\omega-1)W_j-c\binom j2=\omega W_j.
\end{aligned}
\]
The first $jc$ entries of $\mathbf p$ are exactly the first $j$
upper tiers of the $c$ reference cliques, whose sum is $\omega W_j$.
The largest $jc$ entries of $\mathbf q$ have a sum at least that of
$U_j$. Thus
\begin{equation}
\sum_{i=1}^{jc}q_i\geq\sum_{i=1}^{jc}p_i
\qquad (1\leq j\leq r).
\label{EqMajorizationAtClassBoundaries}
\end{equation}

If $t>0$, the partial-class bound gives
\[
|\{e\in E(G):e\cap U\ne\varnothing\}|
=d(rc+t)-e_G(U)
\geq c\left(rd-\binom r2\right)+t(d-r)=W.
\]
All edges incident with $U$ precede the other edges in our ordering,
so all $W$ selected matches have their winners in $U$.
At most $e_G(U)$ selected matches are internal to $U$; consequently,
\[
\begin{aligned}
\sum_{v\in U}q_v
&\geq d(rc+t)+(\omega-1)W-e_G(U)\geq\omega W
=\sum_{i=1}^{rc+t}p_i.
\end{aligned}
\]
This gives the additional prefix inequality at $rc+t$.

We now verify the remaining prefixes. For a decreasing vector,
its interpolated prefix sums form a concave function. The prefix
sums of $\mathbf p$ are linear between consecutive indices among
$0,c,2c,\ldots,rc$, and between $rc$ and $rc+t$ when $t>0$.
Concavity and the established endpoint inequalities therefore give
\[
\sum_{i=1}^{s}q_i\geq\sum_{i=1}^{s}p_i
\qquad(0\leq s\leq rc+t).
\]
Both configurations have the same number $W$ of decisive matches and
therefore the same total number of points:
\begin{equation}
\sum_{i=1}^N q_i=\sum_{i=1}^N p_i
=Nd+(\omega-2)W.
\label{EqEqualTotalPointsMajorization}
\end{equation}
If $t=0$, the remaining reference scores are all $x-1$, so the same
concavity argument proves the inequalities for all later prefixes.
If $t>0$, the remaining reference scores are in
$\{x-1,x-2\}$. Suppose that $s>rc+t$ is the first position at
which the prefix inequality fails. It must then be that $q_s<p_s$.
Since scores are integers, $q_s\leq p_s-1$; and, since $\mathbf q$
is decreasing and all later reference scores are at least $p_s-1$,
we have $q_i\leq p_i$ for every $i\geq s$. The strict prefix deficit
could never be recovered at $N$, contradicting
\eqref{EqEqualTotalPointsMajorization}.

Thus $\mathbf q$ majorises $\mathbf p$. By the convexity of
$z\mapsto z^2$ and equality of their total points,
\[
\mathit{HHI}_\omega^\ast(G)
\geq\mathit{HHI}_\omega(G, \mathbf q)
\geq\mathit{HHI}_\omega(G_{\mathrm{worst}},\mathbf p).
\]
Since $\mathbf p$ is an optimal score vector of
$G_{\mathrm{worst}}$, this proves the first inequality in
\eqref{EqLowerBoundCompleteComponents} for every $\omega\geq2$.
For the second inequality, take any optimal configuration of $K$ and
repeat it on all $c$ components. If its total points and sum of squared
scores are $P_K$ and $S_K$, the resulting configuration on
$G_{\mathrm{worst}}$ has
\[
\mathit{HHI}_\omega
=\frac{cS_K}{(cP_K)^2}
=\frac1c\mathit{HHI}_\omega^\ast(K)
=\frac{d+1}{N}\mathit{HHI}_\omega^\ast(K).
\]
The optimum on $G_{\mathrm{worst}}$ is at least this value, giving the
second inequality in \eqref{EqLowerBoundCompleteComponents}.

\medskip
\noindent\textbf{Step 4: Uniqueness for $\omega\in\{2,3,4,5\}$.}
We may assume $d,c\ge2$; otherwise, any possible $d$-regular graph is of the form $cK_{d+1}$. Write $H=cK_{d+1}$ and choose the optimal
truncated-cascade configuration on $H$ used above, with parameters $r,t$
and $W$ decisive matches.

First, this configuration has at least two fully upgraded vertices. Indeed,
the all-draw configuration (that has the unique worst possible $\mathit{HHI}$ of $\frac{1}{N}$) cannot be optimal. If it had exactly one fully
upgraded vertex, then $r=0,t=1$. Let $F(s)$ denote the HHI of the
configuration with one fully upgraded vertex in each of $s$ cliques and
only draws elsewhere. Then by \eqref{eq:numpoints}-\eqref{eq:numsquares}
\[
 F(s)=\frac{Nd+sA}{d(N+bs)^2},\qquad
 A=(\omega^2-3)d+1,\quad b=\omega-2.
\]
A direct calculation gives
\[
 F(2)>F(1)
 \quad\Longleftrightarrow\quad
 \frac{A}{d}>
 \frac{b(2+3b/N)}{1-2b^2/N^2}.
\]
Since $N=c(d+1)\ge6$ and $b\in\{0,1,2,3\}$, the right-hand side is
at most $0,45/17,54/7,21$, respectively, whereas $A/d$ is greater than
$1,6,13,22$, respectively (using the bound
$A/d\ge b^2+4b+1$). Hence $F(2)>F(1)$, ruling out $r=0,t=1$.
In particular, if $r=0$, then $t\ge2$.

We use the following elementary colouring observation. If $G$ is not a
union of $K_{d+1}$'s, it has a proper equitable $(d+1)$-colouring with
classes $A_0,A_1,\ldots,A_d$ such that
\[
 e_G(A_0,A_1)\le c,
 \qquad N_G(v)\cap A_0=\varnothing\text{ for some }v\in A_1. \tag{$\dagger$}
\]
To see this, start from any equitable colouring. The average number of
edges between two classes is $c$. Any pair with fewer than $c$ edges,
or with exactly $c$ edges but not forming a perfect matching, gives the
claim (interchanging its classes if necessary). Otherwise, all pairs are
joined by perfect matchings. Since $G$ is not a clique union, there is
an induced path $u-v-w$; its three vertices have distinct colours.
Swapping the colours of $v,w$ preserves properness and equitability,
but $u$ now has no neighbour of $v$'s original colour. By the same
averaging argument, the resulting colouring has a pair satisfying
$(\dagger)$.

Fix this colouring, order the remaining classes greedily as in Steps 1--3,
and distinguish two cases.

\emph{If $W>cd$:}  We then have $r\ge 2$ or $r=1 $ and $t\ge 1$; and by~(9), $W-cd\ge d-1$, so the reference has
$p_{c+1}=\omega(d-1)$. In the edge ordering of Step 3, put the $d$
edges incident with $v$ first among the edges led by $A_1$; if a partial
class $B\subseteq A_1$ is used, choose it to contain $v$, which is
possible because $v$ has no neighbour in $A_0$. All $cd$ edges incident
with $A_0$ are selected first. At least $d-1$ of $v$'s edges are
therefore selected, making its score at least $\omega(d-1)+1$.
Meanwhile, all $c$ vertices of $A_0$ score $\omega d$. Thus the sorted
constructed score vector $\mathbf q$ satisfies
$q_{c+1}>p_{c+1}$; by Steps 1--3 it majorizes $\mathbf p$, and strict
convexity of the square yields a strict HHI inequality.

\emph{If $W\le cd$:} The reference either has $(r,t)=(1,0)$ or
$r=0,t\ge2$. Put $h=c$ in the first case and $h=t$ in the second,
so that $W=hd$. Because $A_0$ is independent, it sends exactly $cd$
edges to the $cd$ vertices outside it. By~$(\dagger)$, at least one
outside vertex has no neighbour in $A_0$; hence another has at least
two neighbours there. Choose $T\subseteq A_0$ of size $h\ge2$
containing two such neighbours. Make every vertex of $T$ win all its
matches and declare all remaining matches draws. For $v\notin T$ put
$a_v=|N_G(v)\cap T|$, then $\sum_{v\notin T}a_v=hd$. The total number of wins and points equal those in the reference ($W=hd$). Furthermore, in the reference, each non-winner has 0 or 1 losses. Therefore, for the difference in the sums of squared scores, we have
\[
 S(G)-S(H)=h(\omega d)^2+\sum_{v\notin T}(d-a_v)^2 - ( h(\omega d)^2  +hd(d-1)^2+(N-hd-h)d^2)=\sum_{v\notin T}\bigl(a_v^2-a_v\bigr)>0,
\]
since some $a_v\ge2$. This again gives a strict HHI inequality.

Thus $\operatorname{HHI}^*_{\omega}(G)>
\operatorname{HHI}^*_{\omega}(cK_{d+1})$ for every
$G\not\cong cK_{d+1}$ and $\omega\in\{2,3,4,5\}$.
\end{proof}

\begin{example}[The clique union need not be the unique minimiser]
\label{ExHighOmegaNonunique}
Let $d=2$, $N=6$, and $\omega=6$. The graphs
$G_{\mathrm{worst}}=2K_3$ and $G=C_6$ are nonisomorphic, but
\begin{equation}
\mathit{HHI}_6^\ast(2K_3)
=\mathit{HHI}_6^\ast(C_6)
=\frac{79}{200}.
\label{EqHighOmegaNonunique}
\end{equation}
Indeed, on either graph, let one vertex win both its matches and
make every other match a draw. The resulting score vector is
$(12,2,2,2,1,1)$, giving $\mathit{HHI}_6=158/400=79/200$.

To see that no configuration on either graph has a higher HHI,
let $W\in\{0,\ldots,6\}$ be the number of decisive matches.
By the general identities for total points and squared scores in \eqref{eq:numpoints} and \eqref{eq:numsquares},
we have
\[
P=12+4W,\qquad S=24+68W-\Gamma_6,
\]
where
\[
\Gamma_6
=\sum_{v\in V}\bigl(36w_v\ell_v+25w_v d_v+d_v\ell_v\bigr)
\geq0.
\]
If $W=2$, either a single vertex wins both matches, yielding
$\Gamma_6=2$, or two vertices win one match each, giving
$\Gamma_6\geq50$. Hence $S\leq158$.
If $W=3$, at least one vertex has exactly one win, so
$\Gamma_6\geq25$ and $S\leq203$.
For the remaining values $W\in\{0,1,4,5,6\}$, the bound
$\Gamma_6\geq0$ suffices. Consequently,
\[
\mathit{HHI}_6
\leq
\max\left\{
\frac{158}{20^2},\ \frac{203}{24^2},\
\max_{W\in\{0,1,4,5,6\}}
\frac{24+68W}{(12+4W)^2}
\right\}
=\max\left\{\frac{79}{200},\frac{203}{576},\frac{37}{98}\right\}
=\frac{79}{200}.
\]
Thus equality in the minimum-HHI comparison does not characterise
$G_{\mathrm{worst}}$ when $\omega\geq6$.
\end{example}

\medskip
The lower bounds break down if $d+1$ does not divide $N$, that is, when the teams cannot be partitioned into full round-robin tournament groups, as shown in the next proposition.

\begin{proposition}[The divisibility assumption is necessary]
\label{PropLowerBoundNondivisibleCounterexample}
The conclusion of Theorem~\ref{Theo3} does not hold for arbitrary $N$.
In particular, let $\omega=5$, let $G=K_{3,3}$, and let
$K=K_4$. Then
\begin{equation}
\mathit{HHI}_5^\ast(K_{3,3})
=
\frac{77}{216}
<
\frac{4}{6}\mathit{HHI}_5^\ast(K_4)
=
\frac{158}{441}.
\label{EqNondivisibleCounterexample}
\end{equation}
\end{proposition}
The proof of Proposition~\ref{PropLowerBoundNondivisibleCounterexample} just calculates the $\mathrm{HHI}^\ast$ values of each graph, and is deferred to Appendix~\ref{Appendix_C}.

\medskip
At first sight, it may be unclear why Theorem~\ref{Theo3} is relevant in practice, since we do not know about any tournament where the teams playing in independent round-robin groups are ranked globally in a single table. However, the third-placed teams of all groups are ranked in the recent FIFA World Cup and UEFA European Championships. Furthermore, the competitive balance of a group stage consisting of round-robin groups is usually computed as the average Herfindahl--Hirschman index of the groups in a season (see, for example, \citet{Triguero-RuizAvila-Cano2023}), which is related to considering the group stage as an incomplete round-robin tournament.

\section{Complexity results} \label{Sec4}

We now turn to the computational complexity of maximising the Herfindahl--Hirschman index in an incomplete round-robin tournament. Section~\ref{Sec41} proves the $\mathcal{NP}$-hardness for every fixed integer $\omega \geq 3$, even under strongly balanced schedules. Section~\ref{Sec42} addresses the case of fixed total points, where the reduction is an easier odd-cycle-transversal encoding.

\subsection{Arbitrary win values under balanced schedules} \label{Sec41}

Throughout this section, a win is worth a fixed integer $\omega \geq 3$, a draw is worth one point, and a loss is worth zero points. We prove hardness under two strong restrictions on the schedule.

\begin{figure}[t!]
\centering
\begin{tikzpicture}[scale=1, every node/.style={font=\small}]
  \foreach \ang/\name/\lab/\shade in {
    90/xp1/$x_1^+$/{15},30/xm2/$x_2^-$/{35},-30/xp3/$x_3^+$/{55},
    -90/xm1/$x_1^-$/{15},-150/xp2/$x_2^+$/{35},150/xm3/$x_3^-$/{55}}
    \node[draw,circle,fill=black!\shade,minimum size=7mm,inner sep=1pt] (\name) at (\ang:2) {\lab};
  \draw (xp1)--(xm2)--(xp3)--(xm1)--(xp2)--(xm3)--(xp1)--cycle;
  \foreach \name in {xp1,xm2,xp3,xm1,xp2,xm3}{
    \draw (\name) -- ($(0,0)!1.35!(\name)$);
  }
  \node at (0,-3) {Variable six-cycle};

  \begin{scope}[xshift=6.2cm]
    \node[draw,circle,fill=black!15,minimum size=8mm] (c1) at (90:1.45) {$c_1$};
    \node[draw,circle,fill=black!35,minimum size=8mm] (c2) at (210:1.45) {$c_2$};
    \node[draw,circle,fill=black!55,minimum size=8mm] (c3) at (330:1.45) {$c_3$};
    \draw (c1)--(c2)--(c3)--(c1)--cycle;
    \draw (c1)--($(0,0)!1.55!(c1)$);
    \draw (c2)--($(0,0)!1.55!(c2)$);
    \draw (c3)--($(0,0)!1.55!(c3)$);
    \node at (0,-3) {Clause triangle};
  \end{scope}
\end{tikzpicture}
\caption{Building blocks of the fully balanced core $F_\varphi$ \\ \vspace{0.2cm}
\footnotesize{\emph{Note}: The three shades indicate the classes $V_1,V_2,V_3$; every short external edge is continued to the corresponding occurrence vertex in the other type of block.}}
\label{FigFullyBalancedCore}
\end{figure}
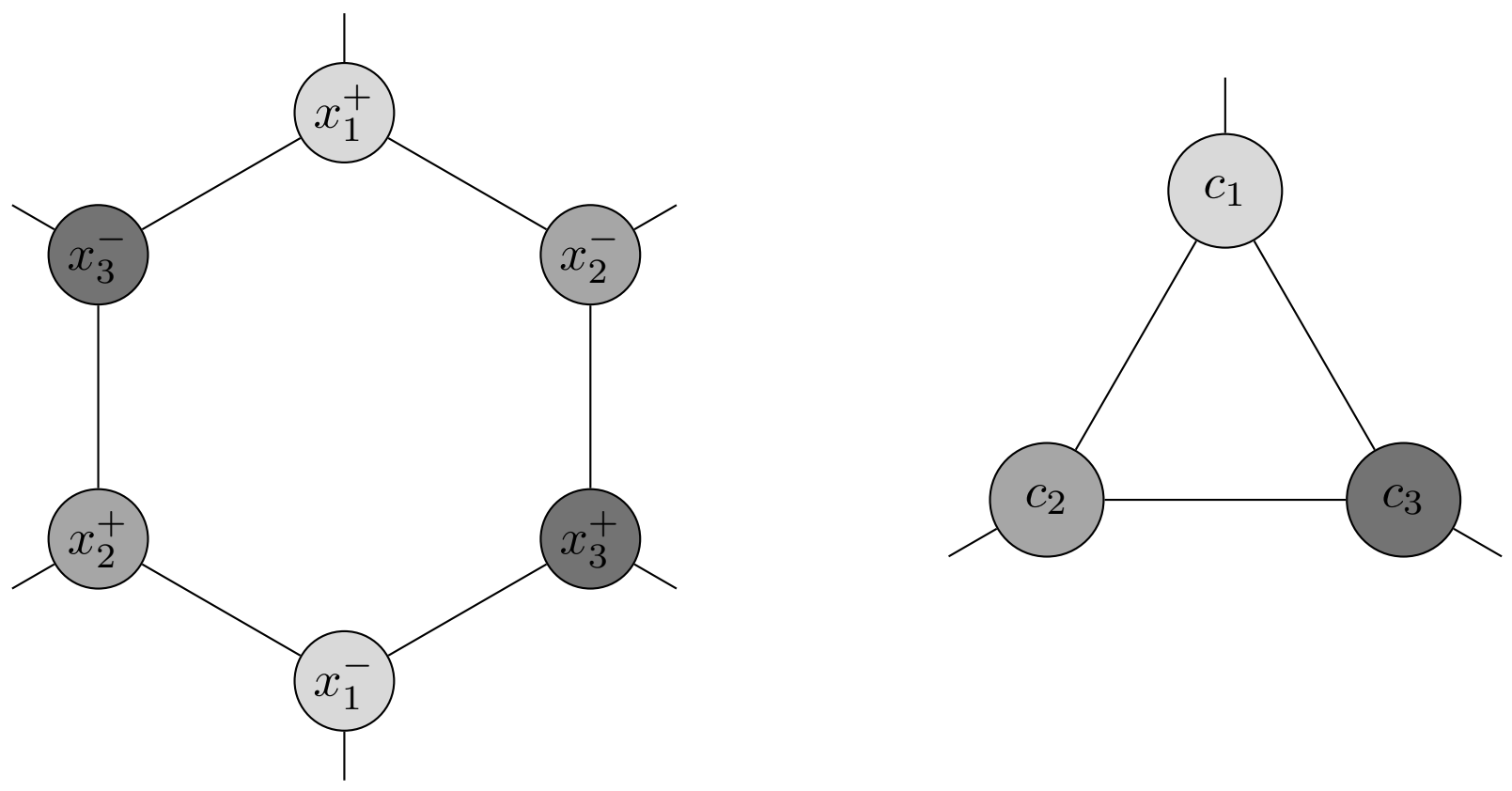

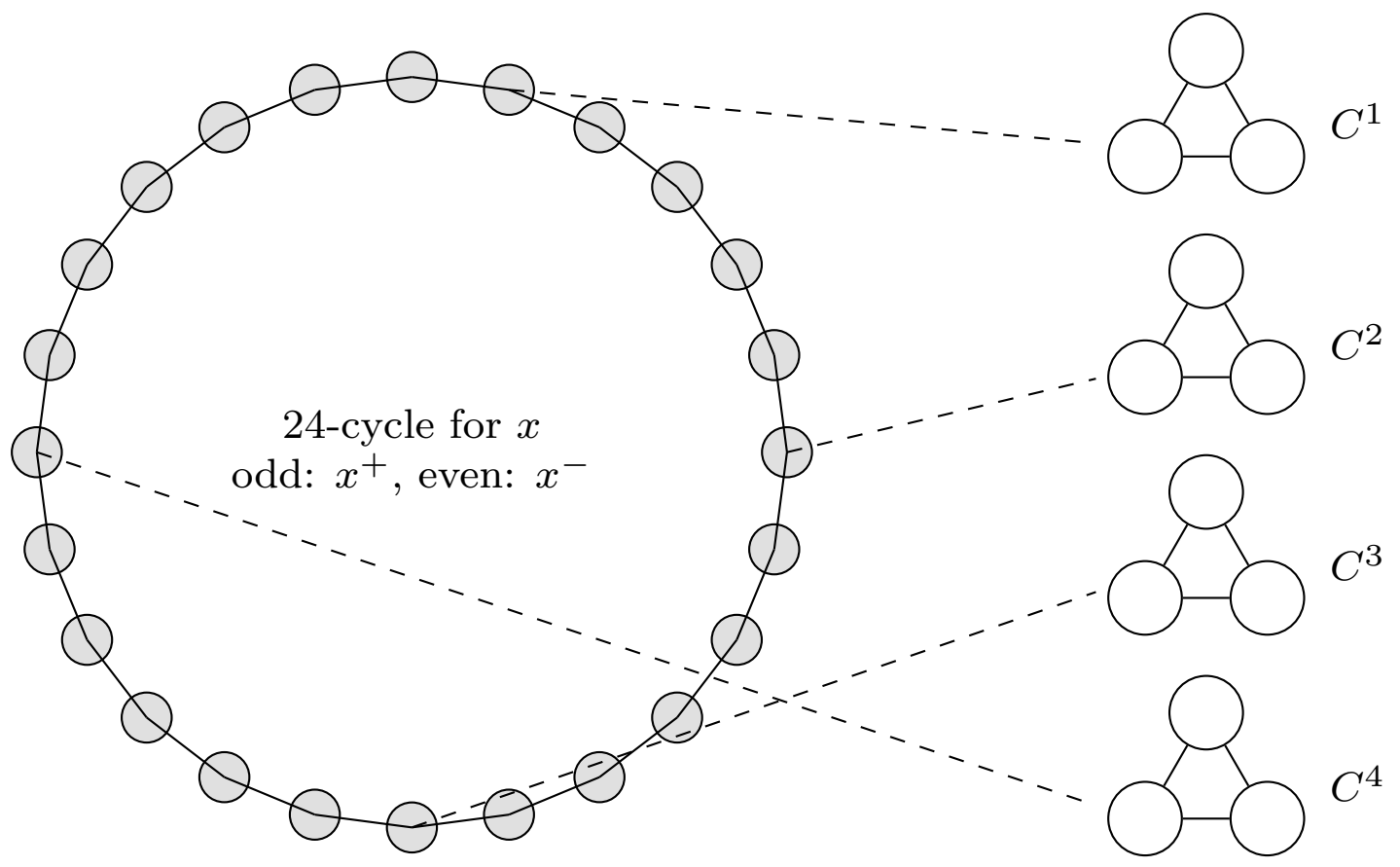
\begin{figure}[t!]
\centering
\begin{tikzpicture}[scale=1, every node/.style={font=\scriptsize}]
  \foreach \i in {1,...,24}{
    \pgfmathsetmacro{\ang}{90-15*(\i-1)}
    \coordinate (v\i) at (\ang:2.55);
    \node[draw,circle,fill=black!12,minimum size=3.4mm,inner sep=0pt] at (v\i) {};
  }
  \foreach \i in {1,...,23}{\pgfmathtruncatemacro{\j}{\i+1}\draw (v\i)--(v\j);}
  \draw (v24)--(v1);
  \node[align=center] at (0,0) {$24$-cycle for $x$\\odd: $x^+$, even: $x^-$};

  \foreach \k/\yy in {1/2.25,2/0.75,3/-0.75,4/-2.25}{
    \begin{scope}[xshift=5.4cm,yshift=\yy cm,scale=.48]
      \node[draw,circle,minimum size=5mm] (a) at (90:1) {};
      \node[draw,circle,minimum size=5mm] (b) at (210:1) {};
      \node[draw,circle,minimum size=5mm] (c) at (330:1) {};
      \draw (a)--(b)--(c)--(a)--cycle;
      \node[draw=none] at (2.15,0) {$C^{\k}$};
    \end{scope}
  }
  \draw[dashed] (v2) -- (4.65,2.1);
  \draw[dashed] (v7) -- (4.65,0.5);
  \draw[dashed] (v13) -- (4.65,-0.95);
  \draw[dashed] (v19) -- (4.65,-2.4);
\end{tikzpicture}
\caption{Schematic view of the externally balanced core $E_\varphi$ \\ \vspace{0.2cm}
}
\label{FigExternallyBalancedCore}
\end{figure}

\noindent
\textbf{Map of the proof:}
The reduction has four layers. First, we rewrite the nonlinear objective as a
linear potential once the total numbers of decisive matches and mixed local
outcomes are recorded. Second, variable cycles and clause triangles encode a
truth assignment: alternating pure winners select one sign of each variable,
while a clause winner must correspond to a true literal. Third, fixed
class-compatible padding gadgets make maximising the $\mathit{HHI}$ equivalent to
maximising this potential on the core. Finally, a complete $k$-fold blow-up
turns parameter one into parameter $k$; and a quantitative gap ensures that
asymmetric configurations of the blow-up cannot erase the separation.
Figures~\ref{FigFullyBalancedCore}--\ref{FigTriangularPrism} illustrate the main building blocks.

\subsubsection{A potential for cubic (3-regular) graphs} \label{Sec411}

Let $\Phi$ be a configuration on a cubic graph. For a vertex $v$, let
$w_v,d_v,\ell_v$ be its numbers of wins, draws, and losses, and let $W$ be
the number of decisive matches. Define
\begin{align}
\Gamma_\omega(\Phi)
&:=
\sum_v\left(
\omega^2w_v\ell_v
+(\omega-1)^2w_vd_v
+d_v\ell_v
\right),
\label{EqSelfContainedGamma}\\
\Psi_q^\omega(\Phi)
&:=qW-\Gamma_\omega(\Phi).
\label{EqSelfContainedPsi}
\end{align}
If the graph has $N$ vertices, then
\begin{align}
P(\Phi)
&=3N+(\omega-2)W,
\label{EqSelfContainedP}\\
S(\Phi)
&=9N+3(\omega^2-2)W-\Gamma_\omega(\Phi).
\label{EqSelfContainedS}
\end{align}
Indeed, $w_v+d_v+\ell_v=3$ gives
\[
w_v^2=3w_v-w_vd_v-w_v\ell_v,
\qquad
d_v^2=3d_v-w_vd_v-d_v\ell_v,
\]
and summing $(\omega w_v+d_v)^2$ yields
\eqref{EqSelfContainedS}. Consequently, for every real number $\lambda$,
\begin{equation}
S(\Phi)-\lambda P(\Phi)
=
3N(3-\lambda)
+
\Psi_{q_\lambda}^\omega(\Phi),
\qquad
q_\lambda:=3(\omega^2-2)-\lambda(\omega-2).
\label{EqSelfContainedLinearisation}
\end{equation}

The following elementary normal form is the only local argument needed in the
section.

\begin{lemma}
\label{LemSelfContainedCanonical}
Let $2\leq q<2(\omega-1)$. On every cubic graph, a configuration maximising
$\Psi_q^\omega$ can be chosen so that every vertex is one of the following:
\[
(3,0,0),\qquad(0,0,3),\qquad(0,2,1),\qquad(0,1,2).
\]
Furthermore, if $\omega\ge 4$, then a configuration maximising $\Psi_1^\omega$ can be chosen such that every vertex is one of the following:
\[
(3,0,0),\qquad(0,0,3),\qquad(0,2,1),\qquad(0,1,2), \qquad (0,3,0).
\]
Thus, every vertex having a win is a pure winner. If $I$ is the set of pure
winners and $R$ is the number of mixed vertices, then
\begin{equation}
\Psi_q^\omega=3q|I|-2R,
\label{EqSelfContainedCanonicalValue}
\end{equation}
and $I$ is an independent set.
\end{lemma}

\begin{proof}
The contribution of a vertex $v$ to $\Psi_q^\omega$ is
\[
\psi_q^\omega(v)
=
qw_v-
\left(
\omega^2w_v\ell_v
+(\omega-1)^2w_vd_v
+d_v\ell_v
\right).
\]

First, let $2\le q<2(\omega -1)$.

\medskip
\noindent
\textbf{Step 1: eliminating vertices with exactly one win.}
Suppose that $v$ has exactly one win, against $u$, and replace this decisive match with a draw. Write the state of $v$ as $(1,d_v,\ell_v)$, where $d_v+\ell_v=2$. At $v$, the change in the potential is
\begin{align*}
\Delta_v
&=
(\omega-1)^2d_v+(\omega^2-1)\ell_v-q\geq
2(\omega-1)^2-q.
\end{align*}
At $u$, one loss is replaced by one draw. The change is
\begin{align*}
\Delta_u
&=(2\omega-1)w_u+d_u-\ell_u+1=(2\omega-2)w_u+4-2\ell_u\geq -2,
\end{align*}
where the last inequality follows from $w_u+d_u+\ell_u=3$. Hence,
\begin{equation}
\Delta_v+\Delta_u
\geq
2\omega^2-4\omega-q.
\label{EqLocalOneWin}
\end{equation}
 If $2\le q<2(\omega-1)$, then it is strictly larger than
\(
2\omega^2-6\omega+2,
\)
which is positive for every $\omega\geq3$. Thus, a maximiser has no vertex with exactly one win.

\medskip
\noindent
\textbf{Step 2: eliminating a vertex with two wins and one draw.}
Suppose that $v$ has state $(2,1,0)$ and draws against $u$. Let $v$ win this match. At $v$, the change is
$\Delta_v=q+2(\omega-1)^2$.

At $u$, one draw is replaced with one loss. Since Step~1 has eliminated vertices with exactly one win, the smallest possible change at $u$ is attained when $u$ has state $(2,1,0)$, and equals $\Delta_u=-4\omega+2$.
Therefore,
\begin{equation}
\Delta_v+\Delta_u
\geq
2\omega^2-8\omega+q+4.
\label{EqLocalTwoWinsDraw}
\end{equation}
This quantity is positive whenever $q>2$, $\omega \ge 3$. For $q=2$, it is nonnegative for every $\omega\geq3$, with equality only when $\omega=3$.

\medskip
\noindent
\textbf{Step 3: eliminating a vertex with two wins and one loss.}
Suppose that $v$ has state $(2,0,1)$. Replace both wins by draws. At $v$, the change is $2(\omega^2-q-1)$.
At either defeated neighbour, replacing one loss by one draw decreases the potential by at most $2$, as calculated in Step~1. Hence, the total change is at least
\begin{equation}
2(\omega^2-q-3).
\label{EqLocalTwoWinsLoss}
\end{equation}
If $2\le q<2(\omega-1)$, then
\(
\omega^2-q-3
>
\omega^2-2\omega-1
>0
\)
for every $\omega\geq3$.

Steps~1--3 show that, when $q>2$, every vertex having a win is a pure winner. The same conclusion holds for $q=2$: among all maximisers, choose one having the largest possible number of pure winners, and apply the non-decreasing operation of Step~2 whenever necessary.

\medskip
\noindent
\textbf{Step 4: eliminating a vertex playing three draws.}
Suppose that $v$ draws all three matches. Since every winner is pure, none of its neighbours has a win. Let $v$ win all three matches. Its contribution increases by $3q$. At each neighbour, replacing one draw by one loss changes the contribution by at least $-2$. Thus, the total change is at least
\begin{equation}
3q-6=3(q-2).
\label{EqLocalPureDraw}
\end{equation}
This is positive for $q>2$. For $q=2$, it is zero; our tie-breaking rule again allows us to perform the operation because it increases the number of pure winners.

Now, assume $\omega\geq4$ and consider the potential $\Psi_1^\omega$. Then, all three modifications in steps 1-3 are
strictly improving, since
\[
2\omega^2-4\omega-1>0,
\qquad
2\omega^2-8\omega+5>0,
\qquad
2(\omega^2-4)>0.
\]
Consequently, in a configuration maximising $\Psi_q^\omega$ or $\Psi_1^\omega$ for $\omega\ge 4$, every winner is
pure. Thus, we also get \eqref{EqSelfContainedCanonicalValue} from \eqref{EqSelfContainedPsi}.
\end{proof}

\begin{corollary}
  If $I$ denotes the set of pure winners, then every edge not incident
with $I$ is a draw, and the contribution of a vertex $v$ is
\begin{equation}
\begin{cases}
3q, & v\in I,\\
-2, & v\notin I\text{ and }|N(v)\cap I|\in\{1,2\},\\
0, & v\notin I\text{ and }|N(v)\cap I|\in\{0,3\}.
\end{cases}
\label{eq:contributions}
\end{equation}
\end{corollary}

\subsubsection{Two balanced satisfiability cores} \label{Sec412}

We reduce from \textsc{Monotone 3-Sat} in which every clause contains three distinct
literals and every variable occurs exactly three times positively and exactly
three times negatively. This problem is $\mathcal{NP}$-complete
\citep{DarmannDocker2021}. Let the formula $\varphi$ have $m$ clauses and
$z$ variables. Counting literal occurrences gives
\begin{equation}
6z=3m,
\qquad
z=\frac{m}{2}.
\label{EqSelfContainedOccurrenceCount}
\end{equation}
Consider the bipartite incidence graph between the clauses and the two signed
literals $x^+,x^-$ of every variable. It is cubic on both sides. By
K\H{o}nig's line-colouring theorem, its edges can be coloured with
$\{1,2,3\}$ so that every clause and every signed literal is incident with
all three colours.

\paragraph{The fully balanced core.}
For every clause, create a triangle, putting the vertex of an occurrence of
colour $j$ into class $V_j$. For every variable $x$, denote its positive and
negative occurrences of colour $j$ by $x_j^+$ and $x_j^-$ and create the
six-cycle
\begin{equation}
x_1^+x_2^-x_3^+x_1^-x_2^+x_3^-x_1^+.
\label{EqSelfContainedSixCycle}
\end{equation}
Put $x_j^+$ and $x_j^-$ into $V_j$ and join every occurrence vertex to its
clause vertex. Call the resulting graph $F_\varphi$.

Every vertex has one neighbour in each of the three classes. Every clause contributes one vertex to each class, and every variable contributes two, so \eqref{EqSelfContainedOccurrenceCount} gives $|V_j|=m+2z=2m$. Thus, $F_\varphi$ is cubic and fully balanced with parameter one.

\begin{lemma}
\label{LemSelfContainedFullyCore}
For $2\leq q<2(\omega-1)$,
\begin{equation}
\max_{\Phi\text{ on }F_\varphi}\Psi_q^\omega(\Phi)
\leq
\left(\frac{15q}{2}-5\right)m.
\label{EqSelfContainedFullyCoreBound}
\end{equation}
Equality holds if and only if $\varphi$ is satisfiable. If $\varphi$ is
unsatisfiable, the right-hand side can be decreased by $3q$.

For $\omega \geq 4$, we also have
\begin{equation}
\max_{\Phi\text{ on }F_\varphi}\Psi_1^\omega(\Phi)
\leq\frac{5m}{2},
\label{EqSelfContainedFullyCoreSmallQ}
\end{equation}
with equality if and only if $\varphi$ is satisfiable; otherwise, the maximum is at most $5m/2-1$.
\end{lemma}

\begin{proof}
For the first assertion, apply Lemma~\ref{LemSelfContainedCanonical}. Every clause triangle contains at most one pure winner and every variable six-cycle at most three, so
\[
|I|\leq m+3z=\frac{5m}{2}.
\]
Every clause triangle contains at least two mixed vertices. If it contains a pure winner, the other two clause vertices lose against it and draw with each other. If it contains no pure winner, all three internal matches are draws; Step~4 implies that every clause vertex must then lose its occurrence match.

Among the $3m$ variable vertices, at most $3z=3m/2$ are pure winners, and at most $m$ are pure losers: every pure-loser occurrence vertex is adjacent to the unique pure winner of its clause triangle. Hence at least $m/2$ variable vertices are mixed and
\[
R\geq2m+\frac{m}{2}=\frac{5m}{2}.
\]
Equation~\eqref{EqSelfContainedCanonicalValue} proves
\eqref{EqSelfContainedFullyCoreBound}.

Equality forces an alternating set of three pure winners in every variable cycle and one pure winner in every clause triangle. Interpret the selected occurrence vertices as false literals. The selected clause vertex is not adjacent to a selected occurrence vertex and therefore represents a true literal. This gives a satisfying assignment. Conversely, a satisfying assignment selects the three false-literal vertices of each variable cycle and one true-literal vertex of every clause triangle, attaining equality. An unsatisfiable formula therefore loses at least one pure winner, and hence at least $3q$ on $\Psi_q^\omega$.

Now let $q=1$ and $\omega\geq4$. By
Lemma~\ref{LemSelfContainedCanonical}, every winner is pure.

We first record a bound that will also be used for the externally balanced
core. Consider a variable cycle of length $2r$, and let $x$ be the number of
pure winners on the cycle. For $j\in\{0,1,2\}$, let $n_j$ be the number of
non-winners having exactly $j$ winning cycle neighbours. Since the winners
form an independent set,
\[
x\leq r,
\qquad
n_1+2n_2=2x.
\]
Before taking the losses towards clause neighbours into account (assuming they are draws), by \eqref{eq:contributions} the contribution of the
cycle is
\[
3x-2(n_1+n_2)=2n_2-x\leq x\leq r.
\]
By \eqref{eq:contributions}, each adjacent pure winner clause vertex can increase this contribution by at
most two, and therefore, if there are $y$ such clause vertices, the cycle
contributes at most $r+2y$. Equality forces $x=n_2=r$, so the winners form
one of the two alternating independent sets.

Applying this estimate with $r=3$, the variable cycles contribute at most
$3z+2Y$, where $Y$ is the total number of pure winner clause vertices. A
clause triangle contributes at most zero if it has no winner and at most
$3-2-2=-1$ if it has a pure winner. Hence,
\[
\Psi_1^\omega\leq 3z+2Y-Y
=\frac{3m}{2}+Y
\leq\frac{5m}{2}.
\]
Equality requires an alternating set of three pure winners in every variable cycle and one pure winner in every clause triangle. As before, interpreting the selected occurrence vertices as false literals and the selected clause vertex as a true literal gives a satisfying assignment, and every satisfying assignment attains equality. Since $\Psi_1^\omega$ is integer-valued, an unsatisfiable instance loses at least one.
\end{proof}

\paragraph{The externally balanced core.}
Let
\[
\sigma_1=(12)(34),
\qquad
\sigma_2=(13)(24),
\qquad
\sigma_3=(14)(23).
\]
For every clause $C$ and every $t\in[4]$, create a triangle $C^t$. If an occurrence has colour $j$, put its vertex in $C^t$ into class $V_{\sigma_j(t)}$ and join it to a variable vertex in class $V_t$.

This already guarantees external balancedness at every clause vertex.

It remains to arrange the variable vertices. A variable vertex is
described by an ordered pair $(a,b)$, where $a$ is its own class and
$b$ is the class of its clause neighbour. Since $a\neq b$, there are
twelve possible ordered pairs. For every variable, we create one
positive and one negative vertex for each such pair, giving $24$
vertices in total.

The availability of these vertices agrees exactly with the clause
construction. For fixed $t$, the three occurrence colours require
the pairs
\[
(t,\sigma_1(t)),\qquad
(t,\sigma_2(t)),\qquad
(t,\sigma_3(t)),
\]
which are precisely the three ordered pairs $(t,b)$ with $b\neq t$.
Thus, as $t$ ranges over $[4]$, every ordered pair $(a,b)$ with
$a\neq b$ is required exactly once for the positive occurrences and
exactly once for the negative occurrences.

Arrange the $24$ variable vertices on a cycle, placing positive
vertices at odd positions and negative vertices at even positions.
Let their own classes be
\begin{equation}
(a_1,\ldots,a_{24})=
(1,4,3,1,2,4,3,2,4,1,2,3,
1,4,2,3,4,1,3,2,1,3,4,2).
\label{EqSelfContainedExternalSequence}
\end{equation}
For every $i\in\mathbb Z_{24}$, define $b_i$ as the unique class not
contained in
$
{a_{i-1},a_i,a_{i+1}}.$
The sequence in
\eqref{EqSelfContainedExternalSequence} has the following two
properties:
\begin{enumerate}
\item

${a_{i-1},a_i,a_{i+1},b_i}=[4]
\qquad
(i\in\mathbb Z_{24});$
\item among the odd positions, and separately among the even
positions, the pairs $(a_i,b_i)$ are exactly the twelve ordered pairs
$(a,b)$ with $a\neq b$ (counted in a cyclic way).
\end{enumerate}
We assign to position $i$ the positive or negative variable vertex
corresponding to the pair $(a_i,b_i)$ and connect it to its clause
neighbour in class $V_{b_i}$.

We can now verify external balancedness at a variable vertex in
position $i$. The vertex itself lies in $V_{a_i}$, its clause
neighbour lies in $V_{b_i}$, and its two cycle neighbours lie in
$V_{a_{i-1}}$ and $V_{a_{i+1}}$. Since
${a_{i-1},a_i,a_{i+1},b_i}=[4],$
these three neighbours lie, one each, in the three classes different
from $V_{a_i}$.

Consequently, every vertex of the resulting graph $E_\varphi$ has
exactly one neighbour in each class different from its own and no
neighbour in its own class. In particular, all four classes are
independent and $E_\varphi$ is cubic.

Finally, every clause contributes three vertices to each class:
class $V_s$ occurs once in each triangle $C^t$ with $t\neq s$.
Moreover, the sequence
\eqref{EqSelfContainedExternalSequence} contains every class six
times, so every variable contributes six vertices to each class.
Thus, all four classes have the same size.

\begin{lemma}
\label{LemSelfContainedExternalCore}
For $2\leq q<2(\omega-1)$,
\begin{equation}
\max_{\Phi\text{ on }E_\varphi}\Psi_q^\omega(\Phi)
\leq(30q-20)m.
\label{EqSelfContainedExternalCoreBound}
\end{equation}
Equality holds if and only if $\varphi$ is satisfiable. In an unsatisfiable
instance, the right-hand side can be decreased by $3q$.

For $\omega\geq4$, we also have
\begin{equation}
\max_{\Phi\text{ on }E_\varphi}\Psi_1^\omega(\Phi)
\leq10m,
\label{EqSelfContainedExternalCoreSmallQ}
\end{equation}
with equality if and only if $\varphi$ is satisfiable; otherwise, the maximum is at most $10m-1$.
\end{lemma}

\begin{proof}
By Lemma~\ref{LemSelfContainedCanonical},
$\Psi_q^\omega=3q|I|-2R$. The $4m$ clause triangles contain at most $4m$
pure winners, and the $z=m/2$ variable cycles contain at most $12z=6m$.
Thus $|I|\leq10m$.

Every clause triangle contains at least two mixed vertices, giving $8m$. Among the $12m$ variable vertices, at most $6m$ are pure winners and at most $4m$ are pure losers, because every pure-loser vertex is adjacent to the unique pure winner of one clause triangle. Therefore, at least $2m$ variable vertices are mixed and $R\geq10m$. This proves the bound \eqref{EqSelfContainedExternalCoreBound}.

Equality forces an alternating independent set in every variable cycle and one pure winner in every clause copy. Odd and even positions represent the positive and negative occurrences, so the alternating set determines a truth value. The selected clause vertex represents a true literal. Conversely, a satisfying assignment attains equality by selecting all false-literal vertices and one true-literal vertex in each clause copy. Unsatisfiability loses at least one pure winner and hence at least $3q$ in $\Psi_q^\omega$.

Now let $q=1$ and $\omega\geq4$. By Lemma~\ref{LemSelfContainedCanonical}, every winner is pure. Apply the cycle estimate from the proof of Lemma~\ref{LemSelfContainedFullyCore} with $r=12$. If $Y$ is the total number of pure winner clause vertices, the
variable cycles contribute at most
\[
12z+2Y=6m+2Y.
\]
As in Lemma~\ref{LemSelfContainedFullyCore}, the clause triangles contribute at most $-Y$. Since $Y\leq4m$,
\[
\Psi_1^\omega\leq6m+Y\leq10m.
\]
Equality forces an alternating independent set in every variable cycle and
one pure winner in every clause copy, and is therefore equivalent to a
satisfying assignment. Since the potential is integer-valued, an
unsatisfiable instance loses at least one.
\end{proof}

\subsubsection{Class-compatible padding} \label{Sec413}

We use the following standard tangent argument.

\begin{lemma}[Tangent padding]
\label{LemSelfContainedTangentPadding}
Let $R$ be a fixed regular graph. Suppose that $R$ has a distinguished configuration whose total number of points and sum of squared scores are $P_0>0$ and $S_0$, respectively. Define \(\lambda:=\frac{2S_0}{P_0}. \) Assume that every configuration of $R$, with corresponding values $(P_R,S_R)$, satisfies \begin{equation} S_R-\lambda P_R \leq S_0-\lambda P_0, \label{EqGeneralTangentCondition} \end{equation} and equality in \eqref{EqGeneralTangentCondition} is possible only if \( (P_R,S_R)=(P_0,S_0). \) Let $H$ be an arbitrary core graph, and add $t$ disjoint copies of $R$ to $H$. For all sufficiently large polynomially bounded $t$, the following statements hold. \begin{enumerate} \item Every configuration maximising the $\mathit{HHI}$ of the resulting disjoint union uses a configuration with values $(P_0,S_0)$ on every copy of $R$. \item Suppose that every copy of $R$ uses such a distinguished configuration. Let $\Phi_1$ and $\Phi_2$ be two configurations of $H$, with values $(P_1,S_1)$ and $(P_2,S_2)$, respectively. If \[ S_1-\lambda P_1>S_2-\lambda P_2, \] then the padded configuration obtained from $\Phi_1$ has strictly larger $\mathit{HHI}$ than the padded configuration obtained from $\Phi_2$. \end{enumerate} The same value of $t$, polynomial in $|V(H)|$, can be chosen simultaneously for all configurations of $H$. \end{lemma}
\medskip
\noindent The proof of Lemma~\ref{LemSelfContainedTangentPadding} is provided in Appendix~\ref{Appendix_D1}.

\begin{figure}[t!]
\centering
\begin{tikzpicture}[scale=1.5, every node/.style={draw,circle,minimum size=7mm,inner sep=1pt}]
  \node[fill=black!10] (a1) at (-1.4,1.0) {$a_1$};
  \node[fill=black!25] (a2) at (0,2.0) {$a_2$};
  \node[fill=black!40] (a3) at (1.4,1.0) {$a_3$};
  \node[fill=black!10] (b1) at (-1.4,-1.0) {$b_1$};
  \node[fill=black!25] (b2) at (0,0) {$b_2$};
  \node[fill=black!40] (b3) at (1.4,-1.0) {$b_3$};
  \draw (a1)--(a2)--(a3)--(a1)--cycle;
  \draw (b1)--(b2)--(b3)--(b1)--cycle;
  \draw[ultra thick] (a1)--(b1);
  \draw[ultra thick] (a2)--(b2);
  \draw[ultra thick] (a3)--(b3);
\end{tikzpicture}
\caption{The triangular prism $\Pi$ and its fully balanced partition \\ \vspace{0.2cm}
\footnotesize{\emph{Note}: The three thick (vertical) matching edges form the three classes.}}
\label{FigTriangularPrism}
\end{figure}
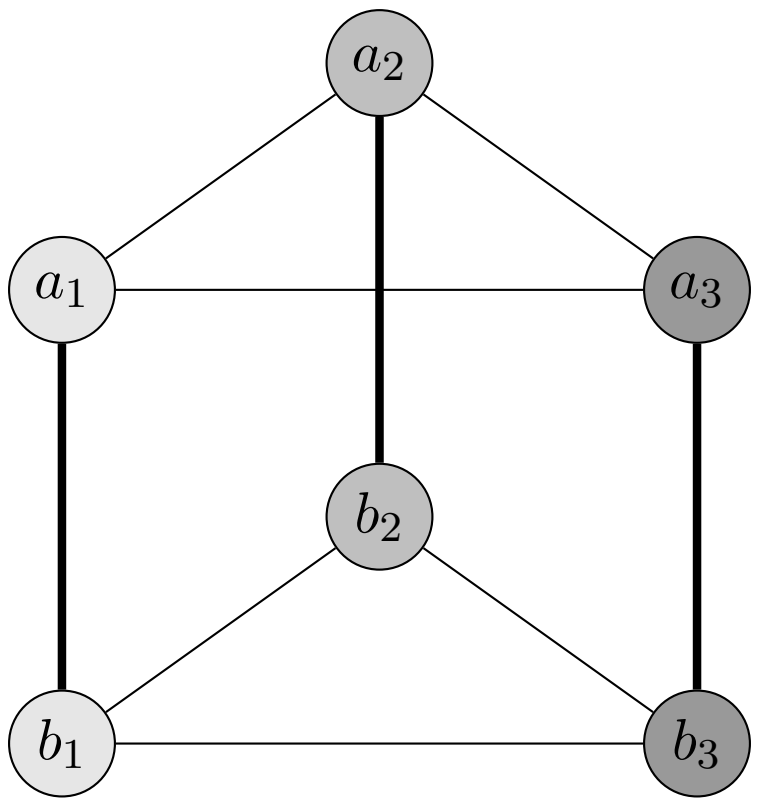

\paragraph{Padding gadgets construction.}
We shall use the following cubic padding gadgets. Partition both $K_{3,3}$
and the triangular prism $\Pi$ (see Figure~\ref{FigTriangularPrism}) into the three edges of a perfect matching.
Let $Q_3$ be $K_{4,4}$ with a perfect matching removed, with the two
nonadjacent vertices $a_i,b_i$ placed in class $i$. In the externally
balanced construction, each vertex of $K_4$ forms a singleton class.

For $\omega\geq4$, put
\[
A_\omega:=3\omega^2-6\omega+5,
\qquad
B_\omega:=3\omega^2-12\omega+7.
\]
Both numbers are positive. Define the fixed cubic padding graphs
\[
R_{\mathrm F}:=
2A_\omega K_{3,3}\mathbin{\dot\cup}\omega B_\omega\Pi,
\qquad
R_{\mathrm E}:=
A_\omega Q_3\mathbin{\dot\cup}\omega B_\omega K_4.
\]
In the distinguished configuration of $R_{\mathrm F}$, every $K_{3,3}$ is oriented bipartitely, and every prism plays only draws. In the distinguished configuration of $R_{\mathrm E}$, every $Q_3$ is oriented bipartitely, and every $K_4$ plays only draws.

\begin{center}
\resizebox{\textwidth}{!}{%
\begin{tabular}{lclccc} \toprule
Structure & $\omega$ & Padding graph and configuration 
& $P_0$ & $S_0$ & $\lambda$ \\
\midrule
fully balanced & $3$ & $K_{3,3}$, pure bipartite orientation
& $27$ & $243$ & $18$ \\
externally balanced & $3$ & $Q_3$, pure bipartite orientation
& $36$ & $324$ & $18$ \\
fully balanced & $\geq4$ & $R_{\mathrm F}$ as above
& $108\omega(\omega-1)(\omega-2)$
& $54\omega(\omega-1)(3\omega^2-7)$
& $(3\omega^2-7)/(\omega-2)$ \\
externally balanced & $\geq4$ & $R_{\mathrm E}$ as above
& $72\omega(\omega-1)(\omega-2)$
& $36\omega(\omega-1)(3\omega^2-7)$
& $(3\omega^2-7)/(\omega-2)$ \\
\bottomrule
\end{tabular}%
}
\end{center}

We verify the hypotheses of Lemma~\ref{LemSelfContainedTangentPadding}. For
$\omega=3$, the pure bipartite configurations maximise $\Psi_3^\omega$ on
$K_{3,3}$ and $Q_3$, with respective values $27$ and $36$. For
$\omega\geq4$,
Lemma~\ref{LemSelfContainedCanonical} gives
\[
\max\Psi_1^\omega(K_{3,3})=9,
\qquad
\max\Psi_1^\omega(Q_3)=12,
\qquad
\max\Psi_1^\omega(\Pi)=
\max\Psi_1^\omega(K_4)=0.
\]
The first two maxima are attained only by a pure bipartite orientation, up to
reversing the two sides, while the last two are attained only by the all-draw
configuration. Indeed, in $K_{3,3}$ a set of two pure winners leaves three
mixed vertices, while three pure winners form one bipartition side and leave
no mixed vertex. In $Q_3$, an independent set of four vertices must be one
bipartition side. In $\Pi$, one pure winner creates three mixed vertices and
two pure winners create four; in $K_4$, one pure winner creates three mixed
vertices. Thus, equality fixes the displayed pair $(P_0,S_0)$ in every row.

For the last two rows, a direct calculation gives
\[
\frac{2S_0}{P_0}
=\frac{3\omega^2-7}{\omega-2}
\qquad\text{and hence}\qquad
q_\lambda=1.
\]
For the first two rows, $\lambda=18$ and $q_\lambda=3$. Consequently,
\eqref{EqSelfContainedLinearisation} shows that the distinguished
configuration satisfies the strict tangent condition, and
Lemma~\ref{LemSelfContainedTangentPadding} applies.

\begin{theorem}[Parameter one] \label{Theo4}
For every fixed integer $\omega\geq3$, computing the maximum $\mathit{HHI}$ is
$\mathcal{NP}$-hard on each of the following graph classes:
\begin{enumerate}
    \item Cubic graphs with a fully balanced partition into three equal-size
    classes with parameter one;
    \item Cubic graphs with an externally balanced partition into four
    equal-size independent classes with parameter one.
\end{enumerate}
\end{theorem}

\begin{proof}
For
$\omega=3$, use the first or second row of the tangent-padding table and the
$q=3$ part of Lemma~\ref{LemSelfContainedFullyCore} or
Lemma~\ref{LemSelfContainedExternalCore}. For $\omega\geq4$, use the
corresponding $R_{\mathrm F}$ or $R_{\mathrm E}$ row and the $q=1$ part of
the relevant core lemma. Equation~\eqref{EqSelfContainedLinearisation} and
Lemma~\ref{LemSelfContainedTangentPadding} reduce maximising $\mathit{HHI}$ to the
corresponding core potential, which is $\mathcal{NP}$-hard. Every core and
every padding component have equal class sizes and the required
parameter-one partition, so their disjoint union preserves the restriction.
\end{proof}

We need one quantitative observation to pass from parameter one to parameter two. Replacing the satisfiability core by a fixed number $L$ of disjoint copies multiplies its strict potential gap by $L$ and preserves all class conditions. The preceding padding calculations therefore give the following.

\begin{lemma}[Amplified gap]
\label{LemSelfContainedAmplifiedGap}
For every fixed integer $\omega\geq 3$ and every constant $A>0$,
either reduction in Theorem~\ref{Theo4} can be arranged to produce a graph $H$ on $n$ vertices and a threshold $\tau$ such that
\[
\text{\rm YES: }\mathit{HHI}_\omega^*(H)\geq\tau,
\qquad
\text{\rm NO: }\mathit{HHI}_\omega^*(H)
\leq
\tau-\frac{A}{n^2}.
\]
\end{lemma}

\medskip
\noindent The proof of Lemma~\ref{LemSelfContainedAmplifiedGap} is provided in Appendix~\ref{Appendix_D2}.

\begin{definition}[Complete $k$-fold blow-up]
Let $G=(V,E)$ be a graph and let $k\geq1$ be an
integer. The \emph{complete $k$-fold blow-up}
$B_k(G)$ is obtained by replacing each vertex
$v\in V$ with an independent set of $k$ clones
$v^1,\ldots,v^k$, and each edge $uv\in E$ with
all $k^2$ edges between the clones of $u$ and
$v$. Formally,
\[
V(B_k(G))=V\times[k],
\quad 
E(B_k(G))=
\bigl\{\{(u,i),(v,j)\}:uv\in E,\ i,j\in[k]\bigr\}.
\]
In particular, if $G$ is $d$-regular on $N$
vertices, then $B_k(G)$ is $kd$-regular on
$kN$ vertices and has $k^2|E|$ edges.
\end{definition}

\begin{lemma}
\label{LemSelfContainedTwoBlowupApproximation}
Let $H$ be a $d$-regular graph on $n$ vertices and let
$h^*=\mathit{HHI}_\omega^*(H)$. Then
\begin{equation}
\frac{h^*}{k}
\leq
\mathit{HHI}_\omega^*(B_k(H))
\leq
\frac{h^*}{k}
+
\frac{\omega(\omega-2)^2}{4n^2k}.
\label{EqSelfContainedTwoBlowupApproximation}
\end{equation}
\end{lemma}

\begin{proof}
The symmetric lift of an optimal configuration of $H$ has $k^2$ times the
total points and $k^3$ times the sum of squared scores, proving the lower
bound.

For the upper bound, take an arbitrary configuration $\widehat\Phi$ of
$B_k(H)$. Independently choose one clone of every original vertex and project
the selected clone edges to a random configuration $\Phi$ of $H$. Then
\begin{equation}
\mathbb{E}[W(\Phi)]=\frac{\widehat W}{k^2},
\qquad \mathbb{E}[P(\Phi)]=\frac{P(\widehat\Phi)}{k^2}, \qquad
\mathbb{E}[\Gamma_\omega(\Phi)]
\leq
\frac{\widehat\Gamma_\omega}{k^3}.
\label{eq:expectations}
\end{equation}
The third inequality follows by expanding each product
$w_v\ell_v,w_vd_v,d_v\ell_v$: a projected ordered pair of distinct original
edges is sampled with probability $1/k^3$, while the corresponding clone
product may contain additional nonnegative pairs arising from the same
original edge. 

Put $a:=\omega^2-2$. Since $H$ has $n$ vertices and degree $d$,
\begin{align*}
\mathbb{E}[S(\Phi)]
&=
nd^2
+
da\,\mathbb{E}[W(\Phi)]
-
\mathbb{E}[\Gamma_\omega(\Phi)]\geq
nd^2
+
da\frac{\widehat W}{k^2}
-
\frac{\widehat\Gamma_\omega}{k^3}.
\end{align*}
On the other hand, $B_k(H)$ has $nk$ vertices and degree $dk$, so
\begin{align*}
\frac{S(\widehat\Phi)}{k^3}
&=
\frac{
(nk)(dk)^2+(dk)a\widehat W-\widehat\Gamma_\omega
}{k^3}=
nd^2
+
da\frac{\widehat W}{k^2}
-
\frac{\widehat\Gamma_\omega}{k^3}.
\end{align*}
Taking expectations,
\begin{equation}
\frac{S(\widehat\Phi)}{k^3}
\leq
h^\ast\mathbb{E}[P(\Phi)^2].
\label{EqProjectedHHIStart}
\end{equation}
From \eqref{EqProjectedHHIStart} and
\(
\mathbb{E}[P(\Phi)^2]
=
\mathbb{E}[P(\Phi)]^2+\operatorname{Var}(P(\Phi)),
\)
we obtain
\begin{align*}
\frac{S(\widehat\Phi)}{k^3}
&\leq
h^\ast
\left[
\frac{P(\widehat\Phi)^2}{k^4}
+
\operatorname{Var}(P(\Phi))
\right].
\end{align*}
Since $S(\Phi)\leq h^*P(\Phi)^2$ for every projection,
\[
\frac{S(\widehat\Phi)}{P(\widehat\Phi)^2}
\leq
\frac{h^*}{k}
+
\frac{k^3h^*\operatorname{Var}(P(\Phi))}{P(\widehat\Phi)^2}.
\]
Now $P=dn+(\omega-2)W$, so
\[
\operatorname{Var}(P)=(\omega-2)^2\operatorname{Var}(W).
\]
Writing $W$ as the sum of the decisiveness indicators of the $nd/2$ edges,
vertex-disjoint indicators are independent, while every variance and every absolute covariance is at most $1/4$. Since there are $n\binom d2$ adjacent edge pairs,
\[
\operatorname{Var}(W)
\leq
\frac{nd}{8}+\frac{nd(d-1)}{4}
\leq
\frac{nd^2}{4}.
\]
Furthermore, $h^*\leq\omega/n$, because
$S\leq(\omega d)P$ and $P\geq nd$, and
$P(\widehat\Phi)\geq k^2nd$. Substitution gives the stated error term.
\end{proof}

\begin{theorem}[Parameter $r$]
\label{Theo5}
For every fixed integer $\omega\geq3$ and $r\ge 1$, computing the maximum $\mathit{HHI}$ is
$\mathcal{NP}$-hard on each of the following graph classes:
\begin{enumerate}
    \item $3r$-regular graphs with a fully balanced partition into three
    equal-size classes with parameter $r$;
    \item $3r$-regular graphs with an externally balanced partition into four
    equal-size independent classes with parameter $r$.
\end{enumerate}
\end{theorem}

\begin{proof}
Use Lemma~\ref{LemSelfContainedAmplifiedGap} with
$A>\omega(\omega-2)^2/4r$ and construct $B_r(H)$. Put all clones of a vertex into its original class. Every original neighbour in a class creates $r$ clone neighbours in that class, so the parameter becomes $r$; independent classes remain independent. In a YES-instance, the symmetric lift has $\mathit{HHI}$ at least $\tau/r$. In a NO-instance,
\[
\mathit{HHI}_\omega^*(B_r(H)) \leq \frac{1}{r}\left(\tau-\frac{A}{n^2}\right) + \frac{\omega(\omega-2)^2}{4n^2r} < \frac{\tau}{r}.
\]
Thus $\tau/r$ separates the two cases.
\end{proof}


\subsection{Fixed total points} \label{Sec42}

We now move to the case when the total number of points is fixed: either draws are prohibited and a win is worth any fixed $\omega>0$ points, or draws are allowed under the $(2,1,0)$ rule. In the latter case, an optimum contains no draws. Indeed, if a draw $uv$ gives the teams $x+1$ and $y+1$ points with $x\geq y$ from their other matches, replacing it by a win of $u$ changes the sum of squared scores by
\[
(x+2)^2+y^2-(x+1)^2-(y+1)^2=2(x-y+1)>0.
\]
Thus, in both settings we may orient every edge and maximise
\[
A(\Phi):=\sum_{v\in V}w_v^2,
\qquad
\mathit{HHI}_\omega(\Phi)=\frac{A(\Phi)}{|E(G)|^2}.
\]

\begin{theorem}[Fixed total points, parameter one] \label{Theo6}
In either fixed-total-points setting, computing the maximum $\mathit{HHI}$ is $\mathcal{NP}$-hard on each of the following classes:
\begin{enumerate}
    \item Cubic graphs with a fully balanced partition into three classes
    with parameter one;
    \item Cubic graphs with an externally balanced partition into four
    independent classes with parameter one.
\end{enumerate}
\end{theorem}

\begin{proof}
For a cubic graph on $N$ vertices, let $\ell_v$ be the indegree of $v$. Since
$w_v+\ell_v=3$ and $\sum_vw_v=3N/2$,
\[
A(\Phi)=\frac{9N}{2}-\sum_vw_v\ell_v.
\]
A vertex contributes two to the last sum exactly when it has both an incoming
and an outgoing edge. If $M(\Phi)$ is the set of such mixed vertices, then
\[
A(\Phi)=\frac{9N}{2}-2|M(\Phi)|.
\]
The pure winners and pure losers form two independent sets, so $M(\Phi)$ is
an odd cycle transversal. Conversely, from any odd cycle transversal $X$,
orient the bipartite graph $G-X$ with bipartition $(L,R)$ from $L$ to $R$, orient the edges between $L$ and $X$ towards $X$ and the ones between $X$ and $R$ towards $R$. Finally, orient the
remaining edges within $X$ arbitrarily. Hence
\begin{equation}
A^*(G)=\frac{9N}{2}-2\tau(G),
\label{EqSelfContainedOrientationOCT}
\end{equation}
where $\tau(G)$ is the minimum odd cycle transversal size.

We reduce from Linear \textsc{Nae-3-Sat} in which every clause has three distinct literals and every variable occurs exactly three times positively and three times negatively; this problem is $\mathcal{NP}$-complete \citep{DarmannDockerDorn2024}. Use the two balanced cubic cores $F_\varphi$ and $E_\varphi$ defined in Section~\ref{Sec412}; the construction itself is independent of the scoring system.

Every clause triangle forces at least one vertex into an odd cycle
transversal. If $\varphi$ is NAE-satisfiable, colour each variable cycle by literal truth values and, in every clause triangle, delete the vertex of one literal while retaining two literals of different truth values. Colour the two surviving clause vertices opposite to their occurrence neighbours. The remaining graph is bipartite. Conversely, a transversal meeting the lower bound deletes exactly one vertex from every clause triangle and no variable vertex; each intact alternating variable cycle determines a truth value, and the two surviving adjacent clause vertices force two literals of the clause
to have different truth values. Therefore,
\[
\tau(F_\varphi)=m
\quad\Longleftrightarrow\quad
\varphi\text{ is NAE-satisfiable},
\]
and, since there are four clause triangles per clause in $E_\varphi$,
\[
\tau(E_\varphi)=4m
\quad\Longleftrightarrow\quad
\varphi\text{ is NAE-satisfiable}.
\]
The first graph has $6m$ vertices and the second has $24m$ vertices. By \eqref{EqSelfContainedOrientationOCT}, their orientation thresholds are $25m$ and $100m$, respectively. Dividing by $|E(G)|^2$ gives the corresponding $\mathit{HHI}$ thresholds.
\end{proof}

We now extend this result to parameter $r$.

\begin{lemma}
\label{LemSelfContainedExactOrientationBlowup}
For every regular graph $G$,
$A^*(B_k(G))=k^3A^*(G)$.
\end{lemma}

\begin{proof}
Lifting an orientation of $G$ symmetrically gives the lower bound. For the
reverse inequality, suppose $G$ is $d$-regular and put
$C(\Phi):=\sum_vw_v\ell_v$. Since $w_v+\ell_v=d$,
\[
A(\Phi)=d|E(G)|-C(\Phi).
\]
Take any orientation $\widehat\Phi$ of $B_k(G)$, independently choose one
clone of every original vertex, and project to an orientation $\Phi$ of $G$.
Expanding $w_v\ell_v$ into ordered pairs of distinct incident edges gives
\[
\mathbb E[C(\Phi)]\leq\frac{C(\widehat\Phi)}{k^3};
\]
the clone product may additionally count two clone edges corresponding to
the same original edge. Hence
$C(\widehat\Phi)\geq k^3\min_\Phi C(\Phi)$. Since $B_k(G)$ has degree $kd$ and
$k^2|E(G)|$ edges,

\[
A(\widehat\Phi)
\leq k^3d|E(G)|-k^3\min_\Phi C(\Phi)
=k^3A^*(G).
\]
\end{proof}

\begin{theorem}[Fixed total points, parameter $r$] \label{Theo7}
In either fixed-total-points setting, computing the maximum $\mathit{HHI}$ is $\mathcal{NP}$-hard on each of the following classes:
\begin{enumerate}
    \item $3r$-regular graphs with a fully balanced partition into three
    classes with parameter $r\ge 1$;
    \item $3r$-regular graphs with an externally balanced partition into four
    independent classes with parameter $r\ge 1$.
\end{enumerate}
\end{theorem}

\begin{proof}
Apply the complete $r$-fold blow-up to either graph in
Theorem~\ref{Theo6} and put all clones of
a vertex into its original class. The class parameter gets multiplied by $r$ and independent
classes remain independent. By
Lemma~\ref{LemSelfContainedExactOrientationBlowup}, the orientation threshold
is multiplied by $r^3$, while the number of edges is multiplied by $r^2$;
thus, the $\mathit{HHI}$ threshold is divided by $r$.
\end{proof}

Finally, the total number of points is also fixed under the (3,2,1,0) scoring system; this case is discussed in Appendix~\ref{Appendix_E}.

\section{Conclusions} \label{Sec5}


Our paper has focused on computing the maximum Herfindahl--Hirschman index in incomplete round-robin tournaments, represented by regular graphs. We have derived sharp lower (Theorem~\ref{Theo3}) and upper (Theorem~\ref{Theo1}) bounds for it. The former is reached when the tournament consists of round-robin groups, while the latter is reached -- under the $(2,1,0)$ and $(3,1,0)$ scoring systems -- when the set of matches forms a bipartite graph. These bounds are important because determining the maximum $\mathit{HHI}$ is $\mathcal{NP}$-hard even if the teams are divided into pots, and each team plays the same number of matches $r$ against teams from any pot. In particular, finding the maximum is computationally difficult for all parameters $r \geq 1$ if the value of a win is at least three points (Theorem~\ref{Theo5}). This is a highly relevant case in practice, because the (3,1,0) scoring system is used in football, and the number of opponents from a given pot is either $r=1$ (UEFA Conference League) or $r=2$ (UEFA Champions League) in some major tournaments. Even though fixing the denominator of $\mathit{HHI}$ by prohibiting draws or using the (2,1,0) scoring system simplifies the problem, computing the maximum $\mathit{HHI}$ remains $\mathcal{NP}$-hard for all $r\ge 1$ (Theorem~\ref{Theo7}). Both Theorems~\ref{Theo5} and \ref{Theo7} remain valid if the teams do not play against teams from their own pot, as in the Champions Hockey League.

Several questions remain open. First, despite the $\mathcal{NP}$-hardness results, it would be interesting to develop approximation algorithms with stronger guarantees, particularly for balanced tournament graphs. 
Practical methods, such as integer programming formulations and local-search heuristics, could also be valuable for computing or approximating the maximum $\mathit{HHI}$ of actual tournament schedules.
A related issue is the range of maximum $\mathit{HHI}$ in a given tournament. For example, in the UEFA Champions League, there are nine teams in each of four pots, and each team faces exactly two opponents from every pot. Hence, the set of matches could not form a bipartite graph, and the lower and upper bounds derived in Section~\ref{Sec3} can be improved. 

Typical tournament structures might be more interesting than extremal cases. For instance, how close is the maximum $\mathit{HHI}$ of a randomly generated balanced tournament to the theoretical upper and lower bounds? Studying its expected value and concentration, especially as the number of teams grows, could reveal whether highly unequal outcomes are mainly a consequence of exceptional schedules or occur frequently even in randomly designed competitions.  Further extensions could consider unequal numbers of matches and repeated encounters between teams.


\section*{Acknowledgements}
\addcontentsline{toc}{section}{Acknowledgements}

The research was supported by the National Research, Development and Innovation Office under Grants Advanced 152220 and FK 145838, and the J\'anos Bolyai Research Scholarship of the Hungarian Academy of Sciences.

\bibliographystyle{apalike} 
\bibliography{All_references}

\clearpage
\appendix

\section*{Appendix}
\addcontentsline{toc}{section}{Appendix}

\section{$\mathit{HHI}_\omega^\ast$ in round-robin tournaments} \label{Appendix_A}

A single round-robin tournament is a special incomplete round-robin tournament, where $d = N-1$ and $G = K_N = K_{d+1}$.
As discussed in the following, the maximum Herfindahl--Hirschman index is known in this case, even if two teams play more than one match against each other.

If $\omega = 2$, that is, two points are awarded for a win and one point for a draw, $\mathit{HHI}_2^\ast$ occurs when one team wins all of its matches, the second team defeats all teams except for the first team, and so on, down to the weakest team that always loses \citep{OwenRyanWeatherston2007}. As shown in \citet[Appendix]{OwenRyanWeatherston2007}:
\[
\mathit{HHI}_2^\ast = \frac{2(2N - 1)}{3N(N-1)}.
\]

However, this upper bound is invalid for sports with $\omega = 3$, such as in (association) football.
Consider a (single) round-robin tournament with $n=4$ teams. The previous logic suggests that the most unequal point distribution is $(9, 6, 3, 0)$, which gives
\[
\mathit{HHI}_3(q=0) = \left( \frac{9}{18} \right)^2 + \left( \frac{6}{18} \right)^2 + \left( \frac{3}{18} \right)^2 + \left( \frac{0}{18} \right)^2 = \frac{1}{4} + \frac{1}{9} + \frac{1}{36} = \frac{7}{18}.
\]
But, if $q$ ($0 \leq q \leq N-1$) teams play only draws against each other, then the total number of points $P$ is reduced, which might increase $\mathit{HHI}_3$. 
Returning to the example, the possible point distribution $(9, 6, 1, 1)$ leads to $\mathit{HHI}_3(q=2) = 81/289 + 36/289 + 2 \cdot 1/289 = 7/17$, while the point distribution $(9, 2, 2, 2)$ results in $\mathit{HHI}_3(q=3) = 9/25 + 3 \cdot 4/225 = 31/75$. Clearly, $\mathit{HHI}_3(q=3) > \mathit{HHI}_3(q=2) > \mathit{HHI}_3(q=0)$.

According to \citet[Proposition~1]{Avila-CanoRuiz-SepulvedaTriguero-Ruiz2021}:
\[
\mathit{HHI}_3(q) = \frac{2q^3 + \left( -6N + 5/2 \right) q^2 + \left( 6N^2 - 5n + 1/2 \right) q + N(N-1)^2}{\left[ -q^2 + (2N - 1) q + 2N(N-1) \right]^2 / 4}.
\]

$\mathit{HHI}_3(q)$ has a unique maximum for a particular $0 \leq q \leq N-1$ \citep[Theorem~1]{Avila-CanoRuiz-SepulvedaTriguero-Ruiz2021}, and the upper bound of $\mathit{HHI}_3$ is given by $\mathit{HHI}_3^\ast = \mathit{HHI}_3(q^\ast)$, where $q^\ast$ is the smallest integer satisfying the following fifth-degree inequality \citep[Theorem~2]{Avila-CanoRuiz-SepulvedaTriguero-Ruiz2021}:
\begin{eqnarray*}
0 & > & \left[ -\frac{1}{2} q^5 + \left( \frac{5}{2}N + \frac{9}{4} \right) q^4 + \left( -8N^2 + 12N - \frac{7}{2} \right) q^3 + \left( 11N^3 + 27N^2 + 18N - \frac{9}{4} \right) q^2 \right. \\
& & + \left. \left( 8N^4 - 6N^3 - 9N^2 + \frac{15}{2} - \frac{1}{2} \right) q - 4N^5 + 12N^4 - 13N^3 + 6N^2 - N \right].
\end{eqnarray*}

It is worth noting that \citet{Avila-CanoRuiz-SepulvedaTriguero-Ruiz2021} only compare win distributions where the first set of teams win all their matches played against the second set of teams, and all matches played by teams from the second set are draws. Consequently, finding the $q$ optimising $\mathit{HHI}_3(q)$ does not necessarily maximise $\mathit{HHI}_3$ globally. According to Lemma~\ref{lem:complete-helper}, there does exist an optimum based on a win distribution of this form for multiple independent round-robin tournaments.


\clearpage
\section{Proof of Step 2 in Lemma~\ref{lem:complete-helper}} \label{Appendix_B}

Let $P_0$ and $S_0$ be the global total points and sum of squared scores,
respectively, when $v$ draws all its $m$ matches in $R$. Define
\[
R_w:=\sum_{i=1}^w r_i
\qquad\text{and}\qquad
R_m:=\sum_{i=1}^m r_i.
\]
If $v$ defeats precisely the first $w$ vertices, then its score changes from
$m$ to
\[
m+(\omega-1)w,
\]
while each defeated vertex loses one point. Consequently,
\begin{equation}
P(w)=P_0+c_\omega w
\label{EqCompleteSaturationP}
\end{equation}
and
\begin{align}
S(w)
={}&
S_0+
\left[m+(\omega-1)w\right]^2-m^2
+\sum_{i=1}^w\left[(r_i-1)^2-r_i^2\right]
\nonumber\\
={}&
S_0+
2m(\omega-1)w
+(\omega-1)^2w^2
-2R_w+w.
\label{EqCompleteSaturationS}
\end{align}

We compare this value with the two boundary configurations $w=0$ and $w=m$.
Put
\[
t:=\frac{w}{m}.
\]
By \eqref{EqCompleteSaturationS},
\begin{align}
&(1-t)S(0)+tS(m)-S(w)
\nonumber\\
&\qquad=
(\omega-1)^2w(m-w)
-
2\left(\frac{w}{m}R_m-R_w\right).
\label{EqCompleteChordGap}
\end{align}

We now bound the term in parentheses. Let $s_i:=r_i-1$ be the score of
vertex $i$ from matches played entirely inside $R\setminus\{v\}$. The added
one point cancels because the coefficients below sum to zero, and hence
\begin{align}
\frac{w}{m}R_m-R_w
&=
\frac{w}{m}\sum_{i=w+1}^m s_i
-
\frac{m-w}{m}\sum_{i=1}^w s_i.
\label{EqCompleteScoreSeparation}
\end{align}
Regard the first $w$ vertices as a bottom set $B$ and the remaining
$m-w$ vertices as a top set $T$. We maximise the right-hand side of
\eqref{EqCompleteScoreSeparation} match by match.

A match inside $T$ has coefficient $w/m$ at both endpoints, so its
contribution is at most $\omega w/m$. A match inside $B$ has coefficient
$-(m-w)/m$ at both endpoints, so its contribution is at most
$-2(m-w)/m$. Finally, a match between $T$ and $B$ contributes the most when
the vertex in $T$ wins, in which case its contribution is $\omega w/m$.
Therefore,
\begin{align}
\frac{w}{m}R_m-R_w
&\leq
\frac{\omega w}{m}
\left[
\binom{m-w}{2}+w(m-w)
\right]
-
\frac{2(m-w)}{m}\binom{w}{2}
\nonumber\\
&=
\frac{w(m-w)}{2m}
\left[
\omega(m+w-1)-2w+2
\right].
\label{EqCompleteScoreSeparationBound}
\end{align}
Combining \eqref{EqCompleteChordGap} and
\eqref{EqCompleteScoreSeparationBound} gives
\begin{equation}
(1-t)S(0)+tS(m)-S(w)
\geq
\kappa_{\omega,m,w}\,w(m-w),
\label{EqCompleteChordGapLower}
\end{equation}
where
\begin{align}
\kappa_{\omega,m,w}
&:=
(\omega-1)^2
-
\frac{\omega(m+w-1)-2w+2}{m}
\nonumber\\
&=
(\omega-1)^2-\omega
-
\frac{c_\omega(w-1)}{m}.
\label{EqCompleteKappa}
\end{align}

Let
\[
M:=\max\left\{
\frac{S(0)}{P(0)^2},
\frac{S(m)}{P(m)^2}
\right\}.
\]
Since both boundary configurations are configurations of the original
$N$-vertex graph and $N\geq m+1$, \eqref{EqCompleteGlobalUpperBound} implies
\begin{equation}
M
\leq
\min\left\{
1,\frac{C_\omega}{m+1}
\right\}.
\label{EqCompleteBoundaryUpper}
\end{equation}
We claim that for $w\notin \{ 0,m\}$
\begin{equation}
\kappa_{\omega,m,w}\geq M c_\omega^2.
\label{EqCompleteKappaDominates}
\end{equation}

First suppose that $\omega=3$. Then $c_\omega=1$, $C_\omega=2$, and
\[
\kappa_{3,m,w}
=
1-\frac{w-1}{m}
=
\frac{m-w+1}{m}
\geq\frac{2}{m}
>
\frac{2}{m+1}
\geq M.
\]

Now suppose that $\omega\geq4$, so $c_\omega\geq2$. Since $\omega = c_\omega +2$ and $w-1\leq m-2$, we get $(\omega-1)^2-\omega=c_\omega^2+c_\omega-1$ and $$\frac{-c_\omega(w-1)}{m}\ge -c_\omega + \frac{2c_\omega}{m}.$$ Thus,
\begin{align}
\kappa_{\omega,m,w}
&\geq
c_\omega^2-1+\frac{2c_\omega}{m}.
\label{EqCompleteKappaCoarse}
\end{align}
If $m\leq2c_\omega$, the right-hand side is at least $c_\omega^2$, and
\eqref{EqCompleteKappaDominates} follows from $M\leq1$.

If $m>2c_\omega$, then
\[
M\leq\frac{C_\omega}{m+1}
<
\frac{C_\omega}{2c_\omega+1}.
\]
Furthermore,
\begin{equation}
C_\omega\leq c_\omega+1.
\label{EqCompleteCUpper}
\end{equation}
Indeed, $c_\omega=\omega-2$, and hence $c_\omega+1=\omega-1$.
Since $\omega\geq4$,
\begin{align*}
C_\omega\leq c_\omega+1
&\iff
\frac{(\omega^2-2)^2}
{4\omega(\omega-1)(\omega-2)}
\leq\omega-1\\
&\iff
(\omega^2-2)^2
\leq
4\omega(\omega-1)^2(\omega-2).
\end{align*}
Substituting $\omega=c_\omega+2$, the difference between the
right-hand side and the left-hand side becomes
\begin{align*}
&4\omega(\omega-1)^2(\omega-2)-(\omega^2-2)^2\\
&\qquad=
4(c_\omega+2)(c_\omega+1)^2c_\omega
-
(c_\omega^2+4c_\omega+2)^2\\
&\qquad=
\left(
4c_\omega^4+16c_\omega^3
+20c_\omega^2+8c_\omega
\right)\\
&\qquad\quad-
\left(
c_\omega^4+8c_\omega^3
+20c_\omega^2+16c_\omega+4
\right)\\
&\qquad=
3c_\omega^4+8c_\omega^3-8c_\omega-4\\
&\qquad=
(3c_\omega^4-4)
+
8c_\omega(c_\omega^2-1)>0,
\end{align*}
where the final inequality follows from $c_\omega\geq2$. Therefore,
$C_\omega<c_\omega+1$, which proves
\eqref{EqCompleteCUpper}.
Since
\[
(c_\omega+1)c_\omega^2
\leq
(2c_\omega+1)(c_\omega^2-1),
\]
where the difference between the right- and left-hand sides is
$c_\omega^3-2c_\omega-1\geq0$, we obtain
\[
Mc_\omega^2
<
\frac{(c_\omega+1)c_\omega^2}{2c_\omega+1}
\leq
c_\omega^2-1
<
\kappa_{\omega,m,w}.
\]
This proves \eqref{EqCompleteKappaDominates}.

By \eqref{EqCompleteSaturationP},
\[
P(w)=(1-t)P(0)+tP(m),
\]
and a direct expansion gives
\begin{equation}
(1-t)P(0)^2+tP(m)^2
=
P(w)^2+c_\omega^2w(m-w).
\label{EqCompletePointChord}
\end{equation}
Using the definition of $M$, followed by
\eqref{EqCompleteChordGapLower},
\eqref{EqCompleteKappaDominates}, and
\eqref{EqCompletePointChord}, we obtain
\begin{align*}
S(w)
&\leq
(1-t)S(0)+tS(m)
-\kappa_{\omega,m,w}w(m-w)\\
&\leq
M\left[(1-t)P(0)^2+tP(m)^2\right]
-Mc_\omega^2w(m-w)\\
&=
MP(w)^2.
\end{align*}
Therefore,
\[
\frac{S(w)}{P(w)^2}
\leq
\max\left\{
\frac{S(0)}{P(0)^2},
\frac{S(m)}{P(m)^2}
\right\}.
\]
Thus, a highest remaining vertex with both wins and draws is dominated by
one of the two boundary configurations: it may be made to draw every match
or to win every match.

\clearpage
\section{Proof of Proposition~\ref{PropLowerBoundNondivisibleCounterexample}} \label{Appendix_C}

\begin{proof}
We first calculate the optimum on $K_4$. By Lemma~\ref{lem:complete-helper}, an optimal configuration is a truncated
cascade. For $\omega=5$, the four possible bottom-tier sizes
$x\in\{1,2,3,4\}$ give the following score vectors:
\begin{center}
\begin{tabular}{c c c}
\toprule
$x$ & Score vector & $\mathit{HHI}_\omega$\\
\midrule
$1$ & $(15,10,5,0)$ & $\displaystyle 7/18$\\[1mm]
$2$ & $(15,10,1,1)$ & $\displaystyle 109/243$\\[1mm]
$3$ & $(15,2,2,2)$ & $\displaystyle 79/147$\\[1mm]
$4$ & $(3,3,3,3)$ & $\displaystyle 1/4$\\
\bottomrule
\end{tabular}
\end{center}
Consequently,
\begin{equation}
\mathit{HHI}_5^\ast(K_4)=\frac{79}{147}.
\label{EqK4OmegaFiveOptimum}
\end{equation}

It remains to determine the optimum on $K_{3,3}$. Let $W$ be the number of
decisive matches. Since $N=6$, $d=3$, and $\omega=5$, the identities for the number of points and sum of squares in \eqref{eq:numpoints} and \eqref{eq:numsquares} become
\begin{align}
P&=18+3W,
\label{EqK33OmegaFiveP}\\
S&=54+69W-\Gamma_5,
\label{EqK33OmegaFiveS}
\end{align}
where
\begin{equation}
\Gamma_5
=
\sum_{v\in V}
\left(
25w_v\ell_v+16w_vd_v+d_v\ell_v
\right).
\label{EqK33OmegaFiveGamma}
\end{equation}

We show that
\begin{equation}
\frac{S}{P^2}\leq\frac{77}{216}.
\label{EqK33DesiredUpper}
\end{equation}
Equivalently, by \eqref{EqK33OmegaFiveP} and
\eqref{EqK33OmegaFiveS}, it is enough to prove
\begin{equation}
\Gamma_5
\geq
54+69W-\frac{77}{216}(18+3W)^2.
\label{EqK33RequiredGamma}
\end{equation}
The right-hand side of \eqref{EqK33RequiredGamma} is nonpositive for
\[
W\in\{0,1,2,7,8,9\},
\]
so these cases follow immediately from $\Gamma_5\geq0$. For the remaining
values, it equals
\begin{center}
\begin{tabular}{c c c c c}
\toprule
$W$ & $3$ & $4$ & $5$ & $6$\\
\midrule
Required lower bound on $\Gamma_5$
& $\displaystyle\frac98$
& $\displaystyle\frac{55}{6}$
& $\displaystyle\frac{259}{24}$
& $6$\\
\bottomrule
\end{tabular}
\end{center}

If a vertex has exactly one or two wins, then its contribution to
$\Gamma_5$ is at least $32$. Indeed, if $w_v=1$, then
$d_v+\ell_v=2$, and
\[
25\ell_v+16d_v+d_v\ell_v\geq32;
\]
if $w_v=2$, then $d_v+\ell_v=1$, and its contribution is either $32$ or
$50$.

Suppose first that $W=3$. If some winning vertex is not a pure winner, then
$\Gamma_5\geq32$. Otherwise, there is exactly one pure winner. Its three
opponents each have state $(0,2,1)$, while all remaining matches are draws.
Thus, these three opponents contribute two each to $\Gamma_5$, and
\[
\Gamma_5\geq6>\frac98.
\]

If $W\in\{4,5\}$, the wins cannot all belong to pure winners because every
pure winner contributes exactly three wins. Hence, some vertex has one or
two wins, and
\[
\Gamma_5\geq32
>
\max\left\{\frac{55}{6},\frac{259}{24}\right\}.
\]

Finally, let $W=6$. If some winning vertex is not pure, then again
$\Gamma_5\geq32$. Otherwise, there are exactly two pure winners. Pure
winners form an independent set, so in $K_{3,3}$ they must belong to the
same bipartition class. Each vertex in the opposite class then has state
$(0,1,2)$ and contributes two to $\Gamma_5$. Therefore,
\[
\Gamma_5\geq6.
\]
This proves \eqref{EqK33DesiredUpper} in every case.

Equality is attained when two vertices in the same bipartition class win
all their matches and the third vertex in that class draws all its matches.
The score vector is
\[
(15,15,3,1,1,1),
\]
so
\begin{equation}
\mathit{HHI}_5^\ast(K_{3,3})
=
\frac{15^2+15^2+3^2+1^2+1^2+1^2}{36^2}
=
\frac{77}{216}.
\label{EqK33OmegaFiveOptimum}
\end{equation}
Finally,
\[
\frac{4}{6}\mathit{HHI}_5^\ast(K_4)
=
\frac23\cdot\frac{79}{147}
=
\frac{158}{441},
\]
and
\[
\frac{158}{441}-\frac{77}{216}
=
\frac{19}{10584}>0.
\]
Hence, the claimed lower bound fails when $d+1$ does not divide $N$.
\end{proof}

\clearpage
\section{Technical proofs for the arbitrary-$\omega$ hardness reduction} \label{Appendix_D}

\subsection{Proof of the tangent-padding lemma (Lemma~\ref{LemSelfContainedTangentPadding})} \label{Appendix_D1}

\begin{proof} We refer to the $t$ copies of $R$ as the \emph{padding gadgets}. The role of these gadgets is to make the nonlinear $\mathit{HHI}$ objective behave, on the core graph $H$, like the linear objective \[ S_H-\lambda P_H. \] Since $R$ is fixed, it has only finitely many possible pairs $(P_R,S_R)$. By the strictness assumption in \eqref{EqGeneralTangentCondition}, there exists a constant $\delta>0$, depending only on $R$, such that every non-distinguished configuration of $R$ satisfies \begin{equation} (S_0-\lambda P_0)-(S_R-\lambda P_R)\geq\delta. \label{EqTangentDeficit} \end{equation} We prove the two statements separately. \paragraph{Step 1: every padding gadget is distinguished.} Fix an arbitrary configuration of the core $H$, and denote its total points and sum of squared scores by $P_H$ and $S_H$. Suppose that the $t$ padding gadgets together have total points and sum of squared scores \[ P_{\mathrm{pad}} \quad\text{and}\quad S_{\mathrm{pad}}. \] Write \[ P_{\mathrm{pad}}=tP_0+x, \qquad S_{\mathrm{pad}}=tS_0+y. \] Thus, $x$ and $y$ measure the change from using the distinguished configuration on every padding gadget. Let $b$ be the number of padding gadgets that are not distinguished. Summing \eqref{EqTangentDeficit} over the $t$ copies gives \[ \begin{aligned} e &:= t(S_0-\lambda P_0) - (S_{\mathrm{pad}}-\lambda P_{\mathrm{pad}})\\ &= \lambda x-y \geq b\delta. \end{aligned} \] Equivalently, \begin{equation} y=\lambda x-e, \qquad e\geq b\delta. \label{EqYLambdaXE} \end{equation} Because $R$ is fixed, there is also a constant $C>0$ such that changing one padding gadget from its distinguished configuration changes its number of points by at most $C$. Therefore, \begin{equation} |x|\leq Cb. \label{EqXBound} \end{equation} We compare the given configuration with the configuration having the same core configuration but using the distinguished configuration on all $t$ padding gadgets. Their $\mathit{HHI}$ values are, respectively, \[ \frac{tS_0+y+S_H}{(tP_0+x+P_H)^2} \qquad\text{and}\qquad \frac{tS_0+S_H}{(tP_0+P_H)^2}. \] Put \[ A:=tP_0+P_H. \] After cross-multiplication, the distinguished padding is better exactly when 
\begin{align*} D &:= (tS_0+S_H)(A+x)^2 - (tS_0+y+S_H)A^2 >0. \end{align*} Expanding the first square gives \[ \begin{aligned} D &= (tS_0+S_H)(2Ax+x^2)-yA^2. \end{aligned} \] Using \eqref{EqYLambdaXE}, we obtain \[ \begin{aligned} D &= (tS_0+S_H)(2Ax+x^2) -(\lambda x-e)A^2\\ &= eA^2 + xA\bigl(2(tS_0+S_H)-\lambda A\bigr) + x^2(tS_0+S_H). \end{aligned} \] The definition $\lambda=2S_0/P_0$ implies \[ 2tS_0-\lambda tP_0=0. \] Consequently, \[ 2(tS_0+S_H)-\lambda A = 2S_H-\lambda P_H, \] and therefore \begin{equation} D = eA^2 + xA(2S_H-\lambda P_H) + x^2(tS_0+S_H). \label{EqPaddingComparisonExact} \end{equation} The last term in \eqref{EqPaddingComparisonExact} is nonnegative. Using $e\geq b\delta$ and $|x|\leq Cb$, we get \[ \begin{aligned} D &\geq b\delta A^2 - CbA\,|2S_H-\lambda P_H|\\ &= bA \left( \delta A-C|2S_H-\lambda P_H| \right). \end{aligned} \] The quantities $P_H$ and $S_H$ are polynomially bounded in $|V(H)|$. Hence, by choosing $t$ polynomially large, we can guarantee \[ \delta(tP_0+P_H) > C|2S_H-\lambda P_H| \] for every possible configuration of $H$. It follows that $D>0$ whenever $b\geq1$. Thus, if even one padding gadget is non-distinguished, replacing all padding gadgets by distinguished ones strictly increases the $\mathit{HHI}$ while leaving the core unchanged. Therefore, every $\mathit{HHI}$-maximising configuration uses $(P_0,S_0)$ on every padding gadget. \paragraph{Step 2: ordering the core configurations.} Assume now that all $t$ padding gadgets are distinguished. A core configuration with values $(P_i,S_i)$ gives the padded $\mathit{HHI}$ \[ F_t(P_i,S_i) := \frac{tS_0+S_i}{(tP_0+P_i)^2}. \] Consider two core configurations with values $(P_1,S_1)$ and $(P_2,S_2)$. Since all denominators are positive, \[ F_t(P_1,S_1)>F_t(P_2,S_2) \] is equivalent to \begin{align} N_t &:= (tS_0+S_1)(tP_0+P_2)^2 - (tS_0+S_2)(tP_0+P_1)^2 >0. \label{EqCoreComparisonCross} \end{align} Expanding \eqref{EqCoreComparisonCross} according to the powers of $t$ gives \begin{align} N_t ={}& t^2P_0^2 \left[ (S_1-\lambda P_1)-(S_2-\lambda P_2) \right] \nonumber\\ &+ t\left[ S_0(P_2^2-P_1^2) + 2P_0(S_1P_2-S_2P_1) \right] \nonumber\\ &+ \left( S_1P_2^2-S_2P_1^2 \right). \label{EqCoreComparisonExpansion} \end{align} The coefficient of $t^2$ is positive precisely when \[ S_1-\lambda P_1 > S_2-\lambda P_2. \] All remaining terms in \eqref{EqCoreComparisonExpansion} grow at most linearly with $t$. Therefore, whenever the above strict inequality holds, the positive quadratic term dominates for all sufficiently large $t$, and hence \[ F_t(P_1,S_1)>F_t(P_2,S_2). \] It remains to verify that one polynomially bounded value of $t$ works for every pair of core configurations. Since \[ \lambda=\frac{2S_0}{P_0} \] and $P_0,S_0$ are fixed integers, every value $S_H-\lambda P_H$ is a rational number whose denominator divides $P_0$. Therefore, two distinct such values differ by at least $1/P_0$. In particular, whenever the coefficient of $t^2$ in \eqref{EqCoreComparisonExpansion} is nonzero, its absolute value is at least $P_0$. On the other hand, all quantities $P_i$ and $S_i$ are polynomially bounded in $|V(H)|$. Hence, the coefficients of $t$ and $t^0$ in \eqref{EqCoreComparisonExpansion} are also polynomially bounded. A polynomially large value of $t$ consequently makes the quadratic term dominate simultaneously for every pair of core configurations. This proves both statements. \end{proof}

\subsection{Proof of the amplified-gap lemma (Lemma~\ref{LemSelfContainedAmplifiedGap})} \label{Appendix_D2}

\begin{proof}
Apply the construction in Theorem~\ref{Theo4} to the conjunction of
$L$ variable-disjoint renamed copies of $\varphi$, where $L$ is a
sufficiently large constant depending only on $A$ and $\omega$. Let $N$
denote the number of vertices in the resulting unpadded core, and let
$(P_0,S_0)$ be the distinguished pair of the corresponding padding graph
from the table above. Put
\[
\lambda:=\frac{2S_0}{P_0}.
\]
For one copy of the formula, let $\Delta_\omega>0$ denote the strict gap in
the relevant core potential: we may take $\Delta_3=9$ and
$\Delta_\omega=1$ for $\omega\geq4$. By
\eqref{EqSelfContainedLinearisation}, a canonical YES-core configuration
$(P_1,S_1)$ and any NO-core configuration $(P_2,S_2)$ therefore satisfy
\[
(S_1-\lambda P_1)-(S_2-\lambda P_2)
\geq L\Delta_\omega.
\]

Add $t$ copies of the padding graph. The difference between the $\mathit{HHI}$ of the
padded YES-configuration and that of the padded NO-configuration is
\[
\frac{D_t}
{(tP_0+P_1)^2(tP_0+P_2)^2},
\]
where
\begin{align*}
D_t
={}&
 t^2P_0^2
\left[
(S_1-\lambda P_1)-(S_2-\lambda P_2)
\right]\\
&+
 t\left[
S_0(P_2^2-P_1^2)
+
2P_0(S_1P_2-S_2P_1)
\right]\\
&+
S_1P_2^2-S_2P_1^2.
\end{align*}
Since the core is cubic and $\omega$ is fixed, $P_i,S_i=O_\omega(N)$.
Hence, the last two lines have absolute value
$O_\omega(tN^2+N^3)$, while the first line is at least
\[
L\Delta_\omega P_0^2t^2.
\]
Choose $t=K_\omega N^2$, where $K_\omega$ is sufficiently large that
Lemma~\ref{LemSelfContainedTangentPadding} applies and the leading term
dominates the remaining terms. Then
\[
D_t\geq C_\omega Lt^2
\]
for some constant $C_\omega>0$ independent of $L$ and of the input formula.
Moreover, $tP_0+P_i=\Theta_\omega(t)$, so the $\mathit{HHI}$ gap is at least
\[
\frac{C'_\omega L}{t^2}.
\]
The total number of vertices in the padded graph is $n=\Theta_\omega(t)$;
therefore, the gap is at least
\[
\frac{C''_\omega L}{n^2}
\]
for a constant $C''_\omega>0$.

Let $\tau$ be the $\mathit{HHI}$ of the padded canonical YES-configuration. The estimate
holds uniformly against every NO-core configuration, and
Lemma~\ref{LemSelfContainedTangentPadding} ensures that an optimum uses the
distinguished configuration on every padding copy. Thus every NO-instance
has optimum at most
\[
\tau-\frac{C''_\omega L}{n^2}.
\]
Choosing $L$ so that $C''_\omega L\geq A$ proves the claim. Since $L$ is a
constant independent of the input formula, the reduction remains
polynomial.
\end{proof}

\clearpage
\section{Maximum $\mathit{HHI}$ under the (3,2,1,0) scoring system} \label{Appendix_E}

Draws are followed by penalty shootouts in some incomplete round-robin competitions, such as the Champions Hockey League (Section~\ref{Sec2}) or the Leagues Cup (Section~\ref{Sec31}), with the following scoring system: three points for a regular win, two points for a shootout win, one point for a shootout loss, and zero points for a regular loss. Consequently, the total number of points is three times the number of matches.
Maximising $\mathit{HHI}$ under this scoring system is equivalent to the $(2,1,0)$ case, which we have investigated thoroughly in Section~\ref{Sec42}.

\begin{proposition}
For every graph $G$ with at least one edge,
\[
\mathit{HHI}_{(3,2,1,0)}^\ast(G) = \mathit{HHI}_{2}^\ast(G).
\]
Moreover, every configuration maximising $\mathit{HHI}$ under the $(3,2,1,0)$ scoring system contains only regular wins and losses.
\end{proposition}

\begin{proof}
Every match awards exactly three points, regardless of its outcome. Hence, the total number of points is fixed at $3|E(G)|$, and maximising $\mathit{HHI}$ is equivalent to maximising the sum of squared scores.

Consider a match between two teams whose scores from their other matches are $a$ and $b$. If the first team receives $x$ points, the sum of their squared scores is
\[
g(x)=(a+x)^2+(b+3-x)^2.
\]
Since $g$ is strictly convex, for $x\in\{1,2\}$ we have
$g(x)<\max\{g(0),g(3)\}.$
Therefore, replacing any shootout result by a  regular victory strictly increases $\mathit{HHI}$. Consequently, every optimal configuration contains only regular wins and losses.

Let $w_v$ denote the number of wins of team $v$. Its score is then $3w_v$, and thus
\[
\mathit{HHI}_{(3,2,1,0)}^\ast(G)
=
\max_{\Phi}
\frac{\sum_{v\in V}(3w_v)^2}
{(3|E(G)|)^2}
=
\max_{\Phi}
\frac{\sum_{v\in V}w_v^2}{|E(G)|^2},
\]
where $\Phi$ ranges over all orientations of $G$.
The last expression is precisely
$\mathit{HHI}_{2}^\ast(G)$, since an optimum under the $(2,1,0)$ scoring system also contains no draws.
\end{proof}

Hence, all extremal characterisations and computational hardness results established for the $(2,1,0)$ scoring system extend directly to the $(3,2,1,0)$ scoring system.

\end{document}